\documentclass[reqno]{amsart}
\usepackage{amsmath}
\usepackage{amssymb}
\usepackage{amsthm,amsxtra}
\usepackage{cite}
\usepackage{enumerate}
\usepackage[mathscr]{eucal}
\usepackage{float}
\newtheorem{Def}{Definition}[section]
\newtheorem{Th}{Theorem}[section]
\newtheorem{Ex}{Example}[section]

\newtheorem{Lemma}{Lemma}[section]
\newtheorem{Prop}{Proposition}[section]
\newtheorem{Cor}{Corollary}[section]

\newcommand{\kc}{\mathcal{K}}
\newcommand{\ks}{\mathscr{K}}
\newcommand{\rb}{\mathbb{R}}
\newcommand{\nb}{\mathbb{N}}

\newcommand{\gc}{\mathcal{G}}

\newcommand{\ac}{\mathcal{A}}
\newcommand{\af}{\mathfrak{a}}
\newcommand{\bc}{\mathcal{B}}
\newcommand{\bk}{\mathfrak{B}}
\newcommand{\fb}{\mathfrak{b}}

\newcommand{\fc}{\mathcal{F}}

\newcommand{\mf}{\mathfrak{m}}

\newcommand{\oc}{\mathcal{O}}

\newcommand{\pf}{\mathfrak{p}}

\newcommand{\uc}{\mathcal{U}}

\newcommand{\vc}{\mathcal{V}}
\newcommand{\wc}{\mathcal{W}}

\newcommand{\obs}{\oc_{\bk^s}}
\newcommand{\bs}{{\bk^s}}

\newcommand{\gbs}{\gamma_{\bk^s}}
\newcommand{\Gbs}{\Gamma_{\bk^s}}

\DeclareMathOperator{\Sf}{{\sf S}_{fin}}
\DeclareMathOperator{\SF}{\sf S_f}
\DeclareMathOperator{\Sh}{\sf S_h}
\DeclareMathOperator{\Sid}{\sf S_{id}}
\DeclareMathOperator{\S1}{{\sf S}_1}

\DeclareMathOperator{\G1}{{\sf G}_1}

\DeclareMathOperator{\CD}{{\sf CD}}
\DeclareMathOperator{\CL}{{\sf CL}}
\DeclareMathOperator{\g}{{\sf G}}
\newcommand{\bp}{\begin{proof}}
\newcommand{\ep}{\end{proof}}

\newcommand{\bns}{(\bk^n)^s}
\newcommand{\bn}{\bk^n}

\newcommand{\Split}{\sf Split}
\newcommand{\FU}{\sf FU}
\newcommand{\SFU}{\sf SFU}
\newcommand{\FUf}{\sf FU_{fin}}
\newcommand{\CFT}{\sf CFT}
\newcommand{\CSFT}{\sf CSFT}
\newcommand{\CSFTf}{\sf CSFT_{fin}}
\newcommand{\IdFT}{\sf IdFT}
\newcommand{\CT}{\sf CT}
\newcommand{\CST}{\sf CST}
\newcommand{\WS}{\sf WS}
\newcommand{\PWS}{\sf PWS}
\newcommand{\WMS}{\sf WMS}
\newcommand{\PS}{\sf PS}
\newcommand{\MS}{\sf MS}
\begin{document}
\title[On tightness and topological games in bornology]{Certain observations on tightness and topological games in bornology}
\author[D. Chandra, P. Das  and S. Das]{Debraj Chandra$^{\dag}$, Pratulananda Das$^*$ and Subhankar Das$^*$\ }
\address{\llap{$\dag$\,}Department of Mathematics, University of Gour Banga, Malda-732103, West Bengal, India}
\email{debrajchandra1986@gmail.com}

\address{\llap{*\,}Department of Mathematics, Jadavpur University, Kolkata-700032, West Bengal, India}
\email{pratulananda@yahoo.co.in, subhankarjumh70@gmail.com}

\thanks{This work of the third author is supported by University Grants Commission (UGC), New Delhi-110002, India, UGC-NET Senior Research Fellowship (1183/(CSIR-UGC NET DEC.2017)).}

\begin{abstract}
This article is a continuation of our investigations in the function space $C(X)$ with respect to the topology $\tau^s_\bk$ of strong uniform convergence on $\bk$ in line of (Chandra et al. 2020 \cite{dcpdsd} and Das et al. 2022 \cite{pddcsd-3}) using the idea of strong uniform convergence (Beer and Levi, 2009 \cite{bl}) on a bornology. First we focus on the notion of tightness property of $(C(X),\tau^s_\bk)$ and some of its variations such as supertightness, Id-fan tightness and $T$-tightness. Certain situations are discussed when $C(X)$ is a {\rm k}-space with respect to the topology $\tau^s_\bk$. Next the notions of strong $\bk$-open game and $\gbs$-open game on $X$ are introduced and some of its consequences are investigated. Finally, we consider discretely selective property and related games. On $(C(X),\tau^s_\bk)$ several interactions between topological games related to discretely selective property, the Gruenhage game on $(C(X),\tau^s_\bk)$ and certain games on $X$ are presented.
\end{abstract}
\maketitle
\smallskip
\noindent{\bf\keywordsname{}:} {Bornology, selection principles, open $\bs$-cover, $\gamma_\bs$-cover, tightness property, topology of strong uniform convergence, {\rm k}-space, discretely selective, the Gruenhage game, function space $C(X)$.}

\section{Introduction}
For notations and terminologies we follow \cite{arh,mcnt,engelking,hh}.
We start with the definition of bornology on a set $X$. A bornology $\bk$ on $X$ is defined as a family of subsets of $X$ such that it is closed under taking finite unions, is hereditary under inclusion and also it forms a cover of $X$ \cite{hh}. A base $\bk_0$ for a bornology $\bk$ is a subfamily of $\bk$ such that for $B\in \bk$ there is a $B_0\in \bk_0$ satisfying $B\subseteq B_0$. When all members of a base $\bk_0$ are closed (compact), $\bk$ is said to have a closed (compact) base. The family $\fc$ of all finite subsets of $X$ is a bornology which is the smallest bornology on $X$. On the other hand, the family of all nonempty subsets of $X$ is the largest bornology on $X$. For a metric space $(X,d)$, the family of all nonempty $d$-bounded subsets of $X$ and the family $\ks$ of nonempty subsets of $X$ with compact closure are bornologies on $X$.

In the realm of metric spaces, the notion of strong uniform continuity on a bornology \cite{bl}, introduced by Beer and Levi in 2009, is one of the most important notions in respect of bornological investigations in function spaces. Let $(X,d)$ and $(Y,\rho)$ be two metric spaces. A mapping $f : X\rightarrow Y$ is strongly uniformly continuous on a subset $B$ of $X$ if for each $\varepsilon>0$ there is a $\delta >0$ such that $d(x_1, x_2)< \delta$ and $\{x_1, x_2\}\cap B \neq \emptyset$ imply $\rho( f(x_1),f(x_2))<\varepsilon $.
Let $\bk$ be a bornology on $X$. We say that $f$ is strongly uniformly continuous on $\bk$ if $f$ is strongly uniformly continuous on $B$ for every $B \in\bk$. In the same article, the authors had also developed a new topology on the set of all functions from $X$ to $Y$, which was designated as the topology of strong uniform convergence on a bornology. Let $\bk$ be a bornology with a closed base on $X$. Then the topology of strong uniform convergence $\tau_{\bk}^s$ is determined by a uniformity on $Y^X$, the set of all functions from $X$ into $Y$, with a base consisting of all sets of the form
 \[ [B,\varepsilon]^s=\{(f,g): \exists \delta> 0 \,\, \text{for every} \,\, x\in B^{\delta},\rho(f(x),g(x))<\varepsilon \},\]
 for $B\in \bk, \varepsilon> 0$ \cite{bl}.

When $\bk= \fc$ the topology of strong uniform convergence $\tau_{\bk}^s$ becomes finer than the topology of pointwise convergence $\tau_p$. If $\bk$ has a compact base, then $\tau_p\leq \tau^s_\bk\leq \tau_k$. Also if $\bk=\ks$, then $\tau_k=\tau^s_\bk$ (see \cite{cmh}). Further investigations on this topology were subsequently carried out in \cite{cmh, h, hn}.

In another direction, a study of the theory of selection principles and their relations to game theory had been initiated systematically by Scheepers in \cite{cooc1, cooc2}. Moreover, the study of selection principles and related games in function spaces endowed with a suitable topology is one of the most fascinating research areas where focus had been given on characterizations of classical properties in function space in respect of selection principles using open covers. Such investigations had been done with respect to the point-open and the compact-open topologies respectively. For more information on selection principles and its recent expansion, see \cite{sur1,sur2,tb}.

From the bornological point of view, with respect to the topology $\tau^s_\bk$ of strong uniform convergence on $\bk$, the study of selection principles in function spaces had been initiated in \cite{cmk}. Further advancement in this direction has been carried out in \cite{dcpdsd, dcpdsd-2, pddcsd-3}, where selection principles and some classical properties have been elaborately investigated using the topology $\tau^s_\bk$ of strong uniform convergence on $\bk$. In this paper we mainly focus on some variations on tightness property, {\rm k}-spaces, discretely selective property and related games in the function space $(C(X),\tau^s_\bk)$. The paper is arranged as follows.

In Section $3$, we deal with the tightness property and some of its variations such as supertightness, Id-fan tightness and $T$-tightness in function spaces. It is shown that the tightness and the supertightness properties of $(C(X),\tau^s_\bk)$ are interchangeable. The Id-fan tightness and the $T$-tightness of $(C(X),\tau^s_\bk)$ are characterized in terms of bornological covering properties of $X$. Section $4$ is devoted to study on {\rm k}-spaces. We show that whenever $(C(X),\tau^s_\bk)$ is a {\rm k}-space it is equivalent to a selection principle related to bornological covers of $X$. In Section $5$, we introduce the notions of strong $\bk$-open game and $\gbs$-open game on $X$ and obtain their consonances with other classical games on $X$. Later in this section, we study discretely selective property and associated games. Under certain condition on $\bk$, it is shown that $(C(X),\tau^s_\bk)$ is discretely selective. Several interactions between topological games on $(C(X),\tau^s_\bk)$ related to discretely selective property, the Gruenhage game on $(C(X),\tau^s_\bk)$ and certain games on $X$ are also presented.

\section{Preliminaries}
We denote the set of all positive integers by $\nb$ .
%We first write down two classical selection principles formulated in general form in \cite{cooc1,cooc2}.
Let $\ac$ and $\bc$ are two nonempty classes of subsets of an infinite set $S$.

$\S1(\ac,\mathcal B)$: It denotes that for any sequence $\{A_n:n\in \nb \}$, where $A_n\in \ac$ for each $n$, there is an $a_n\in A_n$ for each $n$ for which the sequence $\{a_n:n\in \nb\}$ belongs to $\mathcal B$ \cite{cooc1}.

$\Sf(\ac,\mathcal B)$: It denotes that for any sequence $\{A_n:n\in \nb \}$, where $A_n\in \ac$ for each $n\in \nb$, there is a finite (possibly empty) subset $B_n$ of $A_n$ for each $n$ for which the collection $\bigcup_{n\in \nb}B_n$ belongs to $\mathcal B$ \cite{cooc1}.

Also for $f\in \nb^\nb$, we consider the following selection principle.

$\SF(\ac,\mathcal B)$: It denotes that for any sequence $\{A_n:n\in \nb\}$, where $A_n\in \ac$ for each $n\in \nb$, there is a finite (possibly empty) subset $B_n$ of $A_n$ with $|B_n|\leq f(n)$ for each $n$ for which the collection $\bigcup_{n\in \nb}B_n$ belongs to $\mathcal B$ \cite{tb-2}. When $f$ is the identity function we use the notation $\Sid(\ac,\bc)$.

$\Split(\ac, \bc)$: It denotes that every element $A$ of $\ac$ can be split into two elements of $\bc$ \cite{cooc1}.

\noindent $\binom{\ac}{\bc}$: For each element $A$ of $\ac$, there is a set $B$ such that $B\subseteq A$ and $B\in \bc$ \cite{tb04}.

Let $X$ be a Tychonoff topological space. We say that a space $X$ is discretely selective if for any sequence $\{U_n:n\in \nb\}$ of nonempty open subsets of $X$ there exists a $x_n\in U_n$ for each $n$ for which the set $\{x_n:n\in \nb\}$ is a closed and discrete  \cite{tkachuk-1}.

We now consider the following games in our subsequent study.

$\G1(\ac,\bc)$: This denotes an infinitely long game played by two players ONE and TWO. In the $n$-th round, $n\in \nb$, ONE chooses a set $A_n$ from $\ac$. In response to that TWO chooses an element $b_n\in  A_n$. We say that TWO wins the play $\{A_1,b_1, \dotsc, A_n, b_n, \dotsc \}$ if $\{b_n :n\in \nb\}\in \bc$. Otherwise ONE wins.

For $x\in X$ the \textit{Gruenhage game} $\g(X,x)$ on $X$ is defined as follows. Suppose that an infinitely long game is played by two players ONE and TWO. In the $n$th round, ONE chooses an open neighbourhood $U_n$ of $x$. In response to that TWO chooses a point $x_n$ from $U_n$. ONE wins the play if $\{x_n:n\in \nb\}$ converges to $x$. Otherwise TWO wins \cite{gg}. For simplicity we denote the Gruenhage game $\g((C(X),\tau^s_\bk),\underline{0})$ on $(C(X),\tau^s_\bk)$ by $\g(\tau^s_\bk)$.

The game $\CD(X)$ on $X$ is defined as follows. Suppose that an infinitely long game is played by two players ONE and TWO. In the $n$th round, ONE chooses a nonempty open subset $U_n$ of $X$. In response to that TWO chooses a point $x_n$ from $U_n$. TWO wins the play if the set $\{x_n:n\in \nb\}$ is closed and discrete. Otherwise ONE wins \cite{tkachuk-2}. For simplicity we denote the game $\CD(C(X),\tau^s_\bk)$ on $(C(X),\tau^s_\bk)$ by $\CD(\tau^s_\bk)$.

For $x\in X$ the game $\CL(X,x)$ on $X$ is defined as follows. Suppose that an infinitely long game is played by two players ONE and TWO. In the $n$th round, ONE chooses a nonempty open subset $U_n$ of $X$. In response to that TWO chooses a point $x_n$ from $U_n$. ONE wins the play if $x\in \overline{\{x_n:n\in \nb\}}$. Otherwise TWO wins \cite{tkachuk-2}. For simplicity we denote the game $\CL((C(X),\tau^s_\bk),\underline{0})$ on $(C(X),\tau^s_\bk)$ by $\CL(\tau^s_\bk)$.

A \textit{Markov strategy} for TWO in the game $\G1(\ac,\bc)$ is a function $\sigma$ satisfying $\sigma(A,n)=a$ for some $a\in A$ for $A\in \ac$ and $n\in \nb$. We say this Markov strategy is winning if whenever ONE chooses $A_n\in \ac$ during each round $n\in \nb$, TWO wins the play by choosing $\sigma(A_n,n)$ during each round $n\in \nb$ \cite{ch}.

A \textit{predetermined strategy} for ONE is a strategy which depends only on the round of the game. If this strategy for ONE is winning then ONE will be able to win a game irrespective of what TWO is playing \cite{ch}.

For simplicity and convenience of notation we denote `winning strategy' by $\WS$, `predetermined winning strategy' by $\PWS$, `winning Markov strategy' by $\WMS$, `predetermined strategy' by $\PS$ and `Markov strategy' by $\MS$.

Let $\uc$ be an open cover of $X$ and $X\not\in \uc$. $\uc$ is called an \textit{$\omega$-cover} \cite{cooc1} (a $k$\textit{-cover} \cite{koc-ggh}) of $X$ if for every finite (respectively, compact) subset $A$ of $X$ there is an $U\in \uc$ such that $A\subseteq U$. The symbols $\Omega$ and $\kc$ denote the collection of all $\omega$-covers and $k$-covers of $X$ respectively. An open cover $\uc=\{U_n:n\in \nb\}$ is called a \textit{$\gamma$-cover} \cite{cooc1} (a \textit{$\gamma_k$-cover} \cite{koc-ggh}) of $X$ if $\uc$ is infinite and for every finite (respectively, compact) subset $A$ of $X$ there exists a $n_0\in \nb$ satisfying $A\subseteq U_n$ for all $n\geq n_0$. The symbols $\Gamma$ and $\Gamma_k$ denote the collection of all $\gamma$-covers and $\gamma_k$-covers of $X$ respectively.

\textit{Tightness $t(X)$} of $X$ is the smallest infinite cardinal $\mf$ such that if $A\subseteq X$ and $x\in \overline{A}$ then there exists a $B\subseteq A$ satisfying $x\in \overline{B}$ and $|B|\leq \mf$. When $t(X)=\omega$, it is called \textit{countable tightness}. For convenience we denote countable tightness by $\CT$. The tightness of $(C(X),\tau^s_\bk)$ is denoted by $t(\tau^s_\bk)$ \cite{hn}. \textit{Lindel\"{o}f number} $l(X)$ of $X$ is the smallest infinite cardinal $\mf$ such that every open cover of $X$ has a subcover with cardinality less than or equal to $\mf$. We denote the Lindel\"{o}f number of $(C(X),\tau^s_\bk)$ by $l(\tau^s_\bk)$. For a bornology $\bk$, $L^s(X,\bk)$ is the smallest cardinal $\mf$ such that every open $\bs$-cover of $X$ has an open $\bs$-subcover with cardinality less than or equal to $\mf$ \cite{hn}. Let $x\in X$. A family of subsets $\ac$ of $X$ is said to be a $\pi$-network at $x$ if every neighbourhood of $x$ contains a member from $\ac$. The \textit{supertightness} $st(x,X)$ of $x$ is the smallest cardinal $\mf$ for which every $\pi$-network $\ac$ at $x$ consisting of finite subsets of $X$ has a subfamily $\gc$ of cardinality less than or equal to $\mf$ which is a $\pi$-network at $x$. The supertightness $st(X)$ of $X$ is defined by $st(X)=\omega\cdot \sup\{st(x,X):x\in X\}$ \cite{mm} (see also \cite{sakai-st}). When $st(X)=\omega$, it is called \textit{countable supertightness}. For simplicity we denote the supertightness of $(C(X),\tau^s_\bk)$ by $st(\tau^s_\bk)$ and countable supertightness by $\CST$. $X$ has \textit{countable fan tightness} (respectively, \textit{countable strong fan tightness}) if for every $x\in X$ and for every sequence $\{A_n:n\in \nb\}$ of subsets of $X$ with $x\in \overline{A_n}$ for each $n$, there is a finite set $F_n\subseteq A_n$ (respectively, $x_n\in A_n$) for each $n$ such that $x\in \overline{\cup\{F_n:n\in \nb\}}$ (respectively, $x\in \overline{\{x_n:n\in \nb\}}$) \cite{arh, sakai}. For convenience we denote this property by $\CFT$ (respectively, $\CSFT$). $X$ has \textit{countable strong fan tightness for finite sets} (in short, $\CSFTf$) if for each point $x\in X$ and for every sequence $\{\ac_n:n\in \nb\}$ of $\pi$-networks at $x$, where $\ac_n$ consists of finite subsets of $X$, there exists a sequence $\{A_n:n\in \nb\}$ such that $A_n\in \ac_n$ for each $n$ and $\{A_n:n\in \nb\}$ forms a $\pi$-network at $x$ \cite{sakai2}. $X$ has \textit{Id-fan tightness} if for every $x\in X$ and for every sequence $\{A_n:n\in \nb\}$, where $A_n$ is a subsets of $X$ with $x\in \overline{A_n}$ for each $n$, there is an $F_n\subseteq A_n$ with $|F_n|\leq n$ for each $n$ satisfying $x\in \overline{\cup\{F_n:n\in \nb\}}$ \cite{gftm}. For convenience we denote this property by $\IdFT$. The following implications immediately follow. $\CSFTf\Rightarrow \CSFT \Rightarrow \IdFT \Rightarrow \CFT$. The
\textit{$T$-tightness} $T(X)$ is the smallest infinite cardinal $\mf$ such that if $\{F_\alpha:\alpha<\kappa\}$ is an increasing family of closed subsets and $\kappa$ is a regular cardinal greater than $\mf$, then $\cup\{F_\alpha:\alpha<\kappa\}$ is closed in $X$ \cite{j} (see also \cite{sakai-t}). $X$ is said to be \textit{Fr\'{e}chet–Urysohn} (in short, $\FU$) if for each  $A\subseteq X$ and each $x\in \overline{A}$, there is a sequence $\{x_n:n\in \nb\}$ in $A$ such that $\{x_n:n\in \nb\}$ converges to $x$ \cite{gn}. $X$ is said to be \textit{strictly Fr\'{e}chet-Urysohn} (in short, $\SFU$) if for every $x\in X$ and for every sequence $\{A_n:n\in \nb\}$ of subsets of $X$ with $x\in \overline{A_n}$ for each $n$, there is an $a_n\in A_n$ for each $n$ such that $\{a_n:n\in \nb\}$ converges to $x$ \cite{gn}. $X$ is said to be \textit{Fr\'{e}chet-Urysohn for finite sets} (in short, $\FUf$) if for each $x\in X$ and each $\pi$-network $\ac$ at $x$ which consists of finite subsets of $X$, there is a subfamily $\{A_n:n\in \nb\}$ of $\ac$ such that $\{A_n:n\in \nb\}$ converges to $x$ \cite{rs} (see also \cite{gs}). $X$ is a \textit{k-space} if the closed subsets of $X$ are precisely those subset $A$ such that for every compact set $K\subseteq X$, $A\cap K$ is closed in $K$ \cite{McCoy}. Every $\FU$ space is a {\rm k}-space. Let $\ac$ be a family of infinite subsets of $\nb$. $P(\ac)$ denotes that there is a subset $P$ of $\nb$ such that for each $A\in \ac$, $P\setminus A$ is finite. The symbol $\pf$ denotes the smallest cardinal number $k$ for which the following statement is false: For each family $\ac$ if any finite subfamily of $\ac$ has infinite intersection and $|\ac|\leq k$ then $P(\ac)$ holds \cite{vaughan}.

Throughout the paper we assume that $X$ is an infinite metric space. We now recall some classes of bornological covers of $X$. Let $\bk$ be a bornology with closed base on a metric space $X$. For $B\in\bk$ and $\delta>0$, let $B^{\delta}=\bigcup_{x\in B} S_\delta(x)$, where $S_\delta(x)=\{y\in X:d(x,y)<\delta\}$ and $d$ is a metric on $X$. For every $B\in\bk$ and $\delta>0$, $\overline{B}\in \bk$ and $\overline{B^{\delta}}\subseteq B^{2\delta}$ \cite{cmk}. A cover $\uc$ of $X$ is said to be a \textit{strong $\bk$-cover} (in short, $\bs$\textit{-cover}) of $X$ \cite{cmh} if $X\not\in \uc$ and for every $B\in \bk$ there are $U\in \uc$ and $\delta>0$ satisfying $B^\delta\subseteq U$. When all members of $\uc$ are open sets, $\uc$ is said to be an open $\bs$-cover. The symbol $\oc_{\bk^s}$ denotes the collection of all open $\bs$-covers of $X$. We call $X$  \textit{$\bs$-Lindel\"{o}f} \cite{cmk} if every open $\bs$-cover of $X$ contains a countable $\bs$-subcover. An open cover $\{U_n:n\in \nb\}$ of $X$ is called a $\gamma_{\bk^s}$\textit{-cover} \cite{cmk} (see also \cite{cmh}) if it is infinite and for each $B\in \bk$ there exist a $n_0\in  \nb$ and a sequence $\{\delta_n:n\ge  n_0\}$ of positive real numbers satisfying $B^{\delta_n}\subseteq U_n$ for all $n\ge  n_0$. For a metric space $(X,d)$ and $n\in \nb$ consider the product space $X^n$ endowed with the product metric $d^n$ defined as $d^n((x_1,\dotsc, x_n),(y_1,\dotsc, y_n))=\max\{d(x_1,y_1),\dotsc, d(x_n,y_n)\}$. Let $\bk$ be a bornology with closed base on a metric space $(X,d)$. The product bornology $\bk^n$ on $X^n$ is generated by the collection $\{B^n:B\in \bk\}$ \cite{hh}. For any $B\in \bk$, $\delta>0$ and $n\in \nb$, $(B^\delta)^n=(B^n)^\delta$ \cite{dcpdsd-2}.

We denote the collection $\{U_1\cap \dots \cap U_n:U_i\in \uc_i, i=1,\dots,n\}$ by $\uc_1\wedge \dots \wedge \uc_n$. For $x\in X$ and $A\subseteq X$, we write $d(x,A)=\inf\{d(x,a):a\in A\}$.

Let $X$ and $Y$ be two metric spaces. We denote the set of all continuous functions from $X$ to $Y$ by $C(X, Y )$. The point-open topology $\tau_p$ and the compact-open topology $\tau_k$ are two commonly used topologies on $C(X,Y)$. The symbols $(C(X, Y), \tau_p)$ ($C_p(X)$ when $Y = \rb$) and $(C(X, Y), \tau_k)$ ($C_k(X)$ when $Y = \rb$) denote the respective spaces endowed with the topologies $\tau_p$ and $\tau_k$ on $C(X,Y)$.

A bornology $\bk$ is \textit{local} if $\bk$ contains as a member a neighbourhood of each $x\in X$ \cite{bl, h}.
%Whenever $\bk$ is local, the evaluation map $e:(X,d)\times (C(X,Y),\tau^s_\bk)\rightarrow (Y,\rho)$ is continuous \cite{bl}. It is also known that whenever $\bk$ has a countable closed base, $(C(X,Y),\tau^s_\bk)$ is metrizable \cite{bl}.
Note that $(C(X),\tau^s_\bk)$ is homogeneous hence it is enough to concentrate on the point $\underline{0}$ while studying the local properties of this space. Throughout the paper we consider the following assumptions. If $\bk$ is a bornology on a metric space $X$, then $X\not\in \bk$. Also if $\uc$ is a cover of $X$, then $X\not\in \uc$.

\section{Variations on tightness of $(C(X),\tau^s_\bk)$}
\subsection{Certain observations on tightness, supertightness and Id-fan tightness}
In this section we consider the tightness property of $(C(X),\tau^s_\bk)$ and some of its variations such as the supertightness and the Id-fan tightness.

It is known that for a bornology $\bk$ with closed base on a metric space $X$  $t(\tau^s_\bk)=L^s(X,\bk)$ \cite[Theorem 3.5]{hn} holds. Our first result shows that in the product space $X^n$ endowed with the product bornology $\bn$ with a compact base  $t(\tau^s_\bk)=L^s(X^n,\bk^n)$ holds for all $n\in \nb$. For a similar result for $C_p(X)$ see \cite[Corollary 4.7.3]{mcnt}.
\begin{Th}
\label{Tt-1}
If $\bk$ is a bornology with compact base on a metric space $X$, then $t(\tau^s_\bk)=L^s(X^n,\bk^n)$ for all $n\in \nb$.
\end{Th}
\bp
Let $t(\tau^s_\bk)=\mf$. Let $\uc$ be an open $\bns$-cover of $X^n$. There is an open $\bs$-cover $\vc$ of $X$ such that $\{V^n:V\in \vc\}$ refines $\uc$ \cite[Lemma 5.1]{dcpdsd-2}. For $B\in \bk$ choose a $V_B\in \vc$ and a $\delta>0$ such that $B^{2\delta}\subseteq V_B$. Choose a $f_{B,V_B}\in C(X)$ satisfying $f_{B,V_B}(B^\delta)=\{0\}$ and $f_{B,V_B}(X\setminus V_B)=\{1\}$. Clearly $\underline{0}\in \overline{\{f_{B,V_B}:B\in \bk\}}$. As $t(\tau^s_\bk)=\mf$, there is a $\bk'\subseteq \bk$ with $|\bk'|\leq\mf$ such that $\underline{0}\in \overline{\{f_{B,V_B}:B\in \bk'\}}$. Let $\wc=\{V_B:B\in \bk'\}$. For each $V_B$ choose a $U_B\in \uc$ such that $V^n_B\subseteq U_B$. We claim that $\{U_B:B\in \bk'\}$ is an open $\bns$-cover of $X^n$. Let $B^n\in \bk^n$ for $B\in \bk$. Choose a  $f_{B_1,V_{B_1}}\in [B,1]^s(\underline{0})\cap \{f_{B,V_B}:B\in \bk'\}$. Now $f_{B_1,V_{B_1}}\in [B,1]^s(\underline{0})$ implies $B^\delta\subseteq f^{-1}_{B_1,V_{B_1}}(-1,1)$ for some $\delta>0$ and hence $B^\delta\subseteq V_{B_1}$. Clearly $(B^\delta)^n\subseteq V_{B_1}^n$ and $(B^n)^\delta\subseteq U_{B_1}$. This proves our claim. Hence $L^s(X^n,\bk^n)\leq  t(\tau^s_\bk)$.

Conversely, let $L^s(X^n,\bk^n)=\mf$. Let $A\subseteq C(X)$ with $\underline{0}\in \overline{A}$. By \cite[Lemma 2.2]{cmk}, $\uc_m=\{f^{-1}(-\frac{1}{m},\frac{1}{m}):f\in A\}\setminus \{X\}$ is an open $\bs$-cover of $X$. Let $\vc_m=\{f^{-1}(-\frac{1}{m},\frac{1}{m})^n:f\in A\}$. Clearly $\vc_m$ is an open $\bns$-cover of $X^n$. As $L^s(X^n,\bk^n)=\mf$, there is an $A_m\subseteq A$ with $|A_m|\leq \mf$ such that $\{f^{-1}(-\frac{1}{m},\frac{1}{m})^n:f\in A_m\}$ is an open $\bns$-cover. If $A'=\cup_{m\in \nb}A_m$, then it follows that $\underline{0}\in \overline{A'}$ with $|A'|\leq\mf$. Hence $t(\tau^s_\bk)\leq L^s(X^n,\bk^n)$.
\ep

For $C_p(X)$ Sakai showed that $st(C_p(X))=t(C_p(X))$ \cite[Theorem 2.3]{sakai-st} and $l(C_p(X))\geq st(X^n)$ for each $n\in \nb$ \cite[Theorem 2.1]{sakai-st}. We now present similar observations in the bornological setting.

\begin{Th}
\label{Tt-2}
If $\bk$ is a bornology with closed base on a metric space $X$, then $st(\tau^s_\bk)=t(\tau^s_\bk)$.
\end{Th}
\bp
Let $t(\tau^s_\bk)=\mf$. Let $\ac$ be a $\pi$-network at $\underline{0}$ consisting of finite subsets of $C(X)$. For each $n$ and $A\in \ac$, let $U_n(A)=\cap_{f\in A}f^{-1}(-\frac{1}{n},\frac{1}{n})$ and let $F_n=\{f\in C(X): f(X\setminus U_n(A))=\{1\} \text{ for some } A\in \ac\}$. We first show that $\underline{0}\in \overline{F_n}$. Let $[B,\varepsilon]^s(\underline{0})$ be a neighbourhood of $\underline{0}$, where $B\in \bk$ and $\varepsilon>0$. Choose a $n\in \nb$ with $\frac{1}{n}< \varepsilon$. There is an $A\in \ac$ with $A\subseteq [B,\frac{1}{n}]^s(\underline{0})$. Clearly for every $f\in A$ there is a $\delta_f>0$ such that $B^{\delta_f}\subseteq f^{-1}(-\frac{1}{n},\frac{1}{n})$. Choose $\delta=\min\{\delta_f:f\in A\}$. We now have $B^\delta\subseteq U_n(A)$. Choose a $f\in C(X)$ satisfying $f(B^{\frac{\delta}{2}})=\{0\}$ and $f(X\setminus U_n(A))=\{1\}$. Clearly $f\in F_n\cap [B,\varepsilon]^s(\underline{0})$ and hence $\underline{0}\in \overline{F_n}$. As $t(\tau^s_\bk)=\mf$, there is a $G_n\subseteq F_n$ with $|G_n|\leq\mf$ such that $\underline{0}\in \overline{G_n}$. For each $g\in G_n$ there exists an $A_g\in \ac$ satisfying $g(X\setminus U_n(A_g))=\{1\}$. Consider $\ac_n=\{A_g:g\in G_n\}$. Choose $\gc=\cup_{n\in \nb}\ac_n$. We now show that $\gc$ is a $\pi$-network at $\underline{0}$. Let $[B,\varepsilon]^s(\underline{0})$ be a neighbourhood of $\underline{0}$, where $B\in \bk$ and $\varepsilon>0$. Choose a $n\in \nb$ with $\frac{1}{n}<\varepsilon$. As $\underline{0}\in \overline{G_n}$, choose a $g\in [B,\frac{1}{n}]^s(\underline{0})\cap G_n$. Now $B^\delta\subseteq g^{-1}(-\frac{1}{n},\frac{1}{n})$ for some $\delta>0$. Clearly $B^\delta \subseteq U_n(A_g)$ and $B^\delta\subseteq \cap_{f\in A_g}f^{-1}(-\frac{1}{n},\frac{1}{n})$. Subsequently $f\in [B,\frac{1}{n}]^s(\underline{0})$ for all $f\in A_g$ and $A_g\subseteq [B,\varepsilon]^s(\underline{0})$. Therefore $\gc$ is a $\pi$-network at $\underline{0}$ of cardinality less than or equal to $\mf$. Hence $st(\tau^s_\bk)\leq t(\tau^s_\bk)$.

Conversely,  let $st(\tau^s_\bk)=\mf$. Let $A\subseteq C(X)$ with $\underline{0}\in \overline{A}$. Let $\ac=\{\{f\}:f\in A\}$. Clearly $\ac$ is a $\pi$-network at $\underline{0}$. There is a subfamily $\gc$ of $\ac$ with $|\gc|\leq \mf$ such that $\gc$ is a $\pi$-network at $\underline{0}$. Choose $B=\{f:\{f\}\in \gc\}$. Then $\underline{0}\in \overline{B}$ with $|B|\leq \mf$. Hence  $t(\tau^s_\bk)\leq st(\tau^s_\bk)$.
\ep

Note that $(C(X),\tau^s_\bk)$ has $\CT$ if and only if $X$ is $\bs$-Lindel\"{o}f \cite[Theorem 2.1]{cmk}. Thus we have the following.
\begin{Cor}
\label{Ct-1}
Let $\bk$ be a bornology with closed base on a metric space $X$. Then the following statements are equivalent.\\
$(1)$ $(C(X),\tau^s_\bk)$ is $\CST$.\\
$(2)$ $(C(X),\tau^s_\bk)$ is $\CT$.\\
$(3)$ $X$ is $\bs$-Lindel\"{o}f.
\end{Cor}

\begin{Cor}
\label{Ctk-1}
For a Tychonoff space $X$ $st(C_k(X))=t(C_k(X))$ holds.
\end{Cor}

\begin{Th}
\label{Tt-3}
Let $\bk$ be a bornology with closed base on a metric space $X$. The following statements hold.\\
$(1)$ $l(\tau^s_\bk)\geq st(X^n)$ for all $n\in \nb$.\\
$(2)$ $l(\tau^s_\bk)\geq t(X^n)$ for all $n\in \nb$.
\end{Th}
\bp
We prove only $(1)$. Let $l(\tau^s_\bk)=\mf$. Let $(x_1,\dots,x_n)\in X^n$ and $\ac$ be a $\pi$-network at $(x_1,\dots,x_n)$ consisting of finite subsets of $X^n$. We need to show that there is a subfamily $\gc$ of $\ac$ with $|\gc|\leq \mf$ which is a $\pi$-network at $(x_1,\dots,x_n)$. Assume the contrary.

Let $F=\{x_1,\dots,x_n\}$ and $S=\{f\in C(X):f(F)=\{0\}\}$. Since $S$ is closed in $(C(X),\tau^s_\bk)$, $l(S)\leq \mf$. For each $A\in \ac$ choose a finite set $B_A\in \bk$ such that $A=\{(y_1,\dots, y_n):y_1,\dots, y_n\in B_A\}$. We claim that $S\subseteq \cup\{[B_A,1]^s(\underline{0}):A\in \ac\}$. Let $f\in S$. As $f(F)=\{0\}$, $f^{-1}(-1,1)^n$ is a neighbourhood of $(x_1,\dots,x_n)$. There is an $A\in \ac$ such that $A\subseteq f^{-1}(-1,1)^n$ and hence $B_A\subseteq f^{-1}(-1,1)$. As $B_A$ is finite, choose a $\delta>0$ such that $B_A^\delta\subseteq f^{-1}(-1,1)$. Hence $f\in [B_A,1]^s(\underline{0})$. This proves our claim. Since $l(S)=\mf$, there is a subfamily $\gc$ of $\ac$ with $|\gc|\leq\mf$ such that $S\subseteq \cup\{[B_A,1]^s(\underline{0}):A\in \gc\}$. By our assumption, $\gc$ is not a $\pi$-network at $(x_1,\dots, x_n)$. Therefore there is a neighbourhood $U_1\times \cdots \times U_n$ (say) of $(x_1,\dots, x_n)$ such that $A\nsubseteq U_1\times \cdots \times U_n$ for any $A\in \gc$. Let $V=U_1\cap \dots \cap U_n$. Choose an $h\in C(X)$ satisfying $h(F)=\{0\}$ and $h(X\setminus V)=\{1\}$. Clearly $h\in S$. Now $A\nsubseteq U_1\times \cdots \times U_n$ for any $A\in \gc$ implies that $B_A\nsubseteq V$ for any $A\in \gc$. Therefore $h(x)=1$ for some $x\in B_A$. Thus $h\not\in [B_A,1]^s(\underline{0})$ for any $A\in \gc$. This contradicts that $S\subseteq \cup\{[B_A,1]^s(\underline{0}):A\in \gc\}$. Hence our assumption is false.
\ep

\begin{Cor}
\label{Ct-2}
Let $\bk$ be a bornology with closed base on a metric space $X$. The following statements hold.\\
$(1)$ If $(C(X),\tau^s_\bk)$ is Lindel\"{o}f, then $X^n$ has $\CST$ for all $n\in \nb$.\\
$(2)$ If $(C(X),\tau^s_\bk)$ is Lindel\"{o}f, then $X^n$ has $\CT$ for all $n\in \nb$.
\end{Cor}

\begin{Cor}
\label{Ctk-2}
For a Tychonoff space $X$ the following statements holds.\\
$(1)$ $l(C_k(X))\geq st(X^n)$ for all $n\in \nb$.\\
$(2)$ $l(C_k(X))\geq t(X^n)$ for all $n\in \nb$.
\end{Cor}

In \cite[Theorem 3.13]{gftm}, it is shown that the $\IdFT$ of $C_p(X)$ is equivalent to the fact that every finite power of $X$ has the $C^{\prime\prime}$ property.
We now present a characterization of $\IdFT$ of $(C(X),\tau^s_\bk)$ in terms of the selection principle $\S1(\obs,\obs)$ of $X$. We first need the following lemma.

\begin{Lemma}
\label{LId-1}
Let $\bk$ be a bornology with closed base on a metric space $X$. The following statements are equivalent.\\
$(1)$ $X$ satisfies $\S1(\obs,\obs)$.\\
$(2)$ $X$ satisfies $\Sid(\obs,\obs)$.\\
$(3)$ There is an $h\in \nb^\nb$ such that $X$ satisfies $\Sh(\obs,\obs)$.
\end{Lemma}
\bp
$(1)\Rightarrow (2)\Rightarrow (3)$ is follows easily.

$(3)\Rightarrow (1)$. Assume $h\in \nb^\nb$ to be strictly increasing for which $(3)$ holds. Let $\{\uc_n:n\in \nb\}$ be a sequence of open $\bs$-covers of $X$. Define $\vc_1=\uc_1\wedge \dots \wedge \uc_{h(1)}$ and $\vc_n=\uc_{h(n-1)+1}\wedge\dots \wedge \uc_{h(n)}$ for $n>1$. By \cite[Lemma 3.1]{dcpdsd}, $\vc_n$ is an open $\bs$-cover of $X$. Also there is a $\wc_n\subseteq \vc_n$ with $|\wc_n|\leq h(n)$ such that $\cup\{\wc_n:n\in \nb\}$ is an open $\bs$-cover of $X$. Let $\wc_1=\{V_i:i\leq h(1)\}$ and $\wc_n=\{V_{h(n-1)+1}, \dots, V_{h(n)}\}$ for $n>1$. Now for each $i\leq h(1)$ there is a $U_i\in \uc_i$ satisfying $V_i\subseteq U_i$. Also for each $n\in \nb$ with $n>1$ and each $i\leq h(n)-h(n-1)$ there is a $U_{h(n-1)+i}\in \uc_{h(n-1)+i}$ satisfying $V_{h(n-1)+i}\subseteq U_{h(n-1)+i}$. Thus $\{U_i:i\in \nb\}$ is an open $\bs$-cover of $X$ and $(1)$ holds.
\ep

\begin{Th}
\label{TId-2}
Let $\bk$ be a bornology with closed base on a metric space $X$. The following statements are equivalent.\\
$(1)$ $(C(X),\tau^s_\bk)$ has $\CSFT$.\\
$(2)$ $(C(X),\tau^s_\bk)$ has $\CSFTf$.\\
$(3)$ $(C(X),\tau^s_\bk)$ has $\IdFT$.\\
$(4)$ $X$ satisfies $\S1(\obs,\obs)$.
\end{Th}
\bp
The equivalences $(1)\Leftrightarrow (2)\Leftrightarrow (4)$ are from \cite[Theorem 3.4]{pddcsd-3}. The implication $(1)\Rightarrow (3)$ is easily followed.
%and $(4)\Rightarrow (1)$ follows from \cite[Theorem 2.3]{cmk}.

$(3)\Rightarrow (4)$. Let $\{\uc_n:n\in \nb\}$ be a sequence of open $\bs$-covers of $X$. For each $n\in \nb$ and $B\in \bk$ there exist a $\delta>0$ and a $U\in \uc_n$ such that $B^{2\delta}\subseteq U$. Let $\uc_{n,B}=\{U\in \uc_n: B^{2\delta}\subseteq U\}$. Choose a $f_{B,U}\in C(X)$ satisfying $f_{B,U}(B^\delta)=\{0\}$ and $f_{B,U}(X\setminus U)=\{1\}$. Let $A_n=\{f_{B,U}:B\in \bk, U\in \uc_{n,B}\}$. Clearly $\underline{0}\in \overline{A_n}$. By $(3)$, there is a $C_n\subseteq A_n$ with $|C_n|\leq n$ satisfying $\underline{0}\in \overline{\cup\{C_n:n\in \nb\}}$. Choose $\vc_n=\{U\in \uc_n:f_{B,U}\in C_n\}$. Observe that $|\vc_n|\leq n$. Let $B\in \bk$. For the neighbourhood $[B,1]^s(\underline{0})$ of $\underline{0}$, we have $[B,1]^s(\underline{0})\cap C_n\neq \emptyset$ for some $n$ and $B^\delta\subseteq f^{-1}_{B,U}(-1,1)$ for some $f_{B,U}\in C_n$ and $\delta>0$. Clearly $B^\delta\subseteq U$ for some $U\in \vc_n$. Therefore $\cup\{\vc_n:n\in \nb\}$ is an open $\bs$-cover of $X$. By Lemma \ref{LId-1}, $X$ satisfies $\S1(\obs,\obs)$.

\ep

It is known that $C_k(X)$ has $\CSFT$ (respectively, $\CSFTf$) if and only if $X$ satisfies $\S1(\kc,\kc)$ \cite[Theorem 2.2]{kocfs}, \cite[Corollary 3.3]{pddcsd-3}. We now have the following.
\begin{Cor}
\label{Ctk-3}
For a Tychonoff space $X$ the following statements are equivalent.\\
$(1)$ $C_k(X)$ has $\CSFT$.\\
$(2)$ $C_k(X)$ has $\CSFTf$.\\
$(3)$ $C_k(X)$ has $\IdFT$.\\
$(4)$ $X$ satisfies $\S1(\kc,\kc)$.
\end{Cor}

\subsection{On $T$-tightness of $(C(X),\tau^s_\bk)$}
In \cite[Theorem 2.3]{sakai-t}, for a Tychonoff topological space $X$, Sakai characterized the $T$-tightness of $C_p(X)$ in terms of a covering property of $X^n$, $n\in \nb$. In this section we present a characterization of the $T$-tightness of $(C(X),\tau^s_\bk)$ in terms of a bornological covering property of $X^n$ endowed with the product bornology $\bn$, $n\in \nb$. For two infinite cardinals $\af,\fb$ with $\af\leq \fb$ let $[\af,\fb]$ denote the set of all cardinals $\mf$ with $\af\leq \mf\leq \fb$. The symbol $\mf^+$ denotes the successor cardinal of $\mf$. We first introduce the following notion.
\begin{Def}
\label{DT-1}
{\rm A space $X$ is called a} $[\af,\fb]_r$-$\bs$ {\rm space if every open $\bs$-cover of $X$ whose cardinality is a regular cardinal $\mf\in [\af,\fb]$ has an open $\bs$-subcover of cardinality less than $\mf$.}
\end{Def}
The superscript $r$ denotes the restriction to regular cardinals. A space $X$ is called a $[\af,\infty]_r$-$\bs$ space if $X$ is a $[\af,\fb]_r$-$\bs$ space for all $\fb\geq \af$. When $L^s[X,\bk]\leq \mf$, $X$ is $[\mf^+,\infty]_r$-$\bs$ space. In line of the definition of $T(\tau)$ \cite{sakai-t}, we introduce the following property.
\begin{Def}
\label{DT-2}
{\rm Let $\bk$ be a bornology with closed base on a metric space $X$. $X$ satisfies the property $T(\mf,\bk)$ if $\kappa$ is a regular cardinal greater than $\mf$ and $\{\uc_\alpha:\alpha<\kappa\}$ is an increasing family of open sets so that $\cup_{\alpha<\kappa} \uc_\alpha$ is an open $\bs$-cover of $X$, then there exists a $\alpha<\kappa$ such that $\uc_\alpha$ is an open $\bs$-cover of $X$.}
\end{Def}

\begin{Th}
\label{TT-1}
Let $\bk$ be a bornology with a compact base $\bk_0$ on a metric space $X$. The following statements are equivalent.\\
$(1)$ $X$ is a $[\af,\fb]_r$-$\bs$ space.\\
$(2)$ Every decreasing sequence $\{F_\alpha:\alpha<\mf\}$ of nonempty closed subsets of $X$, where $\mf\in [\af,\fb]$ is regular, has nonempty intersection.
\end{Th}
\bp
$(1)\Rightarrow (2)$. Let $\mf\in [\af,\fb]$ be regular. Let $\{F_\alpha:\alpha<\mf\}$ be a decreasing sequence of nonempty closed subsets of $X$. Define $U_\alpha=X\setminus F_\alpha$. Clearly $U_\alpha\neq X$ for all $\alpha<\mf$. We show that $\cap_{\alpha<\mf} F_\alpha\neq\emptyset$. Assume that $\cap_{\alpha<\mf} F_\alpha=\emptyset$. Observe that $\cap_{\alpha<\mf}(X\setminus U_\alpha)=\emptyset$ and $X=\cup_{\alpha<\mf}U_\alpha$. Let $\uc=\{U_\alpha:\alpha<\mf\}$ and choose $B\in \bk_0$. Since $B$ is compact and $B\subseteq\cup_{\alpha<\mf}U_\alpha$ and $U_\alpha$'s are increasing, there is a $U_\alpha$ such that $B\subseteq U_\alpha$. Moreover, there is a $\delta>0$ such that $B^\delta\subseteq U_\alpha$. Hence $\uc$ is an open $\bs$-cover of $X$. By $(1)$, there is a $\tau<\mf$ such that $\{U_\alpha:\alpha<\tau\}$ is an open $\bs$-cover of $X$. As $U_\alpha$'s are increasing and $\cup_{\alpha<\tau}U_\alpha=X$, we have $U_\tau=X$. This is a contradiction as $U_\tau\neq X$. Hence our assumption is false.

$(2)\Rightarrow (1)$. Let $\mf\in [\af,\fb]$ be regular. Let $\uc=\{U_\alpha:\alpha<\mf\}$ be an open $\bs$-cover of $X$ which has no $\bs$-subcover of cardinality less than $\mf$. Let $\uc_\beta=\{U_\alpha:\alpha<\beta\}$. Clearly $\uc_\beta\subset \uc_\tau$ whenever $\beta<\tau$. Define $F_\beta=X\setminus \cup \uc_\beta$. Clearly $\{F_\beta:\beta<\mf\}$ is a decreasing sequence of nonempty closed subsets of $X$. We have $\cup_{\beta<\mf} (\cup\uc_\beta)=X$ and $\cap_{\beta<\mf}(X\setminus \cup\uc_\beta)=\emptyset$. Thus $\cap_{\beta<\mf}F_\beta=\emptyset$. This contradicts $(2)$. Hence there is a $\beta<\mf$ such that $\uc_\beta$ is an open $\bs$-subcover of $X$.
\ep

\begin{Th}
\label{TT-2}
Let $\bk$ be a bornology with compact base on a metric space $X$. The following statements are equivalent.\\
$(1)$ $X$ satisfies $T(\mf,\bk)$.\\
$(2)$ $X^n$ satisfies $T(\mf,\bn)$ for any $n\in \nb$.
\end{Th}
\bp
$(1)\Rightarrow (2)$. Let $\kappa>\mf$ be regular. Let $\{\uc_\alpha:\alpha<\kappa\}$ be an increasing family of open sets so that $\uc=\cup_{\alpha<\kappa} \uc_\alpha$ is an open $\bns$-cover of $X^n$. By \cite[Lemma 5.1]{dcpdsd-2}, there is an open $\bs$-cover $\vc$ of $X$ such that $\{V^n:V\in \vc\}$ refines $\uc$. Let $\vc_\alpha=\{V\in \vc:V^n\subseteq U \text{ for some }U\in \uc_\alpha\}$. Clearly $\vc_\alpha\subseteq \vc_\beta$ whenever $\alpha<\beta$. By $(1)$, there is a $\alpha<\kappa$ such that $\vc_\alpha$ is an open $\bs$-cover of $X$. We now show that $\uc_\alpha$ is an open $\bns$-cover of $X$. Let $B^n\in \bn$ for some $B\in \bk$. Choose a $V\in \vc_\alpha$ and a $\delta>0$ such that $B^\delta\subseteq V$. Clearly $(B^\delta)^n\subseteq V^n$ and so $(B^n)^\delta\subseteq U$ for some $U\in \uc_\alpha$ with $V^n\subseteq U$. Hence $\uc_\alpha$ is an open $\bns$-cover of $X^n$ and $(2)$ holds.

$(2)\Rightarrow (1)$. Let $\uc=\cup_{\alpha<\kappa}\uc_\alpha$ be an open $\bs$-cover of $X$ and $\uc_\alpha\subseteq \uc_\beta$ whenever $\alpha<\beta$. Let $\vc_\alpha=\{U^n:U\in \uc_\alpha\}$. Now $\vc=\cup_{\alpha<\kappa}\vc_\alpha$ is an open $\bns$-cover of $X^n$. By $(2)$, there is a $\alpha<\kappa$ such that $\vc_\alpha$ is an open $\bns$-cover of $X^n$. Thus $\uc_\alpha$ is an open $\bs$-cover of $X$.
\ep

\begin{Th}
\label{TT-3}
Let $\bk$ be a bornology with compact base on a metric space $X$. The following statements are equivalent.\\
$(1)$ $X$ satisfies $T(\mf,\bk)$.\\
$(2)$ $X^n$ is a $[\mf^+,\infty]_r$-$\bns$ space for any $n\in \nb$.
\end{Th}
\bp
$(1)\Rightarrow (2)$. We first show that $X$ is a $[\mf^+,\infty]_r$-$\bs$ space. Let $\kappa>\mf$ be regular and $\{F_\alpha:\alpha<\kappa\}$ be a decreasing sequence of nonempty closed subsets of $X$. Let $U_\alpha=X\setminus F_\alpha$ and $\uc_\beta=\{U_\alpha:\alpha<\beta\}$. Clearly $\uc_\beta$'s are increasing. Assume that $\cap_{\alpha<\kappa}F_\alpha=\emptyset$. Clearly $\cap_{\alpha<\kappa}(X\setminus U_\alpha)=\emptyset$ and $X=\cup_{\alpha<\kappa}U_\alpha$. Let $\uc=\cup_{\beta<\kappa}\uc_\beta$. Now $\uc$ is an open $\bs$-covers of $X$. By $(1)$, there is a $\beta<\kappa$ such that $\uc_\beta$ is an open $\bs$-cover of $X$. Since $U_\alpha$'s are increasing, we have $U_\beta=X$ and $F_\beta=\emptyset$. This contradicts that $F_\beta\neq \emptyset$. Therefore $\cap_{\alpha<\kappa}F_\alpha\neq \emptyset$. By Theorem \ref{TT-1}, $X$ is a $[\mf^+,\infty]_r$-$\bs$ space. Also by Theorem \ref{TT-2}, $X^n$ satisfies $T(\mf,\bn)$ for any $n\in \nb$. Repeating the above procedure, it is easy to see that $X^n$ is a $[\mf^+,\infty]_r$-$\bns$ space for $n>1$.

$(2)\Rightarrow (1)$. Let $\kappa>\mf$ be regular. Let $\uc=\cup_{\alpha<\kappa}\uc_\alpha$ be an open $\bs$-cover of $X$, where $\uc_\alpha\subseteq \uc_\beta$ if $\alpha<\beta$. Fix a $n\in \nb$. Let $\vc_\alpha=\{U^n:U\in \uc_\alpha\}$. The collection $\vc=\cup_{\alpha<\kappa}\vc_\alpha$ is an open $\bns$-cover of $X^n$ such that $\vc_\alpha\subseteq \vc_\beta$ if $\alpha<\beta$. By $(2)$, there is a $\bns$-subcover of $\vc$ of cardinality less than $\kappa$. As $\vc_\alpha\subseteq \vc_\beta$ whenever $\alpha<\beta$, let $\vc_\tau$ be that $\bns$-subcover of $\vc$, where $\tau<\kappa$. Clearly $\uc_\tau$, $\tau<\kappa$, is an open $\bs$-cover of $X$. Hence $(1)$ holds.
\ep

%Now using Theorem \ref{TT-3} and \cite[Theorem 5.2]{dcpdsd} we obtain the following characterization of $T$-tightness of $(C(X),\tau^s_\bk)$. For the convenience of the readers we outline the proof here.

\begin{Th}
\label{TT-4}
Let $\bk$ be a bornology with compact base on a metric space $X$. The following statements are equivalent.\\
$(1)$ $T(C(X),\tau^s_\bk)\leq \mf$.\\
$(2)$ $X^n$ is a $[\mf^+,\infty]_r$-$\bns$ space for any $n\in \nb$.
\end{Th}

\bp
$(1)\Rightarrow (2)$. Let $T(C(X),\tau^s_\bk)\leq \mf$. We show that $X$ satisfies $T(\mf,\bk)$. Let $\kappa>\mf$ be regular and $\cup_{\alpha<\kappa}\uc_\alpha$ be an open $\bk^s$-cover of $X$ such that $\uc_\alpha\subseteq \uc_\beta$ whenever $\alpha<\beta$. For $B\in \bk$ choose a $\delta>0$ and a $U\in \cup_{\alpha<\kappa} \uc_\alpha$ such that $B^{2\delta}\subseteq U$. Choose a $f_{B,U}\in C(X)$ satisfying $f_{B,U}(B^\delta)=\{0\}$ and $f_{B,U}(X\setminus U)=\{1\}$. Consider $A_\alpha=\{f_{B,U}:U\in \uc_\alpha, B\in \bk\}$ for each $\alpha<\kappa$. Clearly $\{\overline{A_\alpha}:\alpha< \kappa\}$ is an increasing sequence of closed subsets of $C(X)$. By $(1)$, the set $\cup_{\alpha<\kappa}\overline{A_\alpha}$ is closed in $(C(X),\tau^s_\bk)$. Using \cite[Lemma 2.2]{cmk} and the fact that $\cup_{\alpha<\kappa}\uc_\alpha$ is an open $\bk^s$-cover of $X$, we obtain $\underline{0}\in \overline{\cup_{\alpha<\kappa} A_\alpha}\setminus \cup_{\alpha<\kappa} A_\alpha$. It is easy to see that $\underline{0}\in \cup_{\alpha<\kappa}\overline{A_\alpha}$. Thus $\underline{0}\in \overline{A_\alpha}$ for some $\alpha<\kappa$. Again by \cite[Lemma 2.2]{cmk}, $\{f_{B,U}^{-1}(-1,1):f_{B,U}\in A_\alpha\}$ is an open $\bs$-cover of $X$. Since $f_{B,U}^{-1}(-1,1)\subseteq U$, $\uc_\alpha$ is an open $\bs$-cover of $X$. Hence $X$ satisfies $T(\mf,\bk)$. The result now follows from Theorem \ref{TT-3}.

$(2)\Rightarrow (1)$. Let $\kappa>\mf$ be regular and let $\{A_\alpha:\alpha<\kappa\}$ be a sequence of closed subsets of $C(X)$ which are increasing. We show that $A=\cup_{\alpha<\kappa}A_\alpha$ is closed in $C(X)$. For this let $\underline{0}\in \overline{A}$. Consider the collection $\uc_{\alpha,n}=\{f^{-1}(-\frac{1}{n},\frac{1}{n}):f\in A_\alpha\}$ and $\uc_n=\cup_{\alpha<\kappa}\uc_{\alpha,n}$ for each $n\in \nb$. Clearly $\uc_n$ is an open $\bs$-cover of $X$. By Theorem \ref{TT-3}, $X$ satisfies $T(\mf,\bk)$. Hence there exists a $\alpha_n<\kappa$ for each $n$ such that $\uc_{{\alpha_n},n}$ is an open $\bs$-cover. Let $\beta=\sup\{\alpha_n:n\in \nb\}$.
Now it can be easily checked that $\uc_{\beta,n}$ is an open $\bk^s$-cover of $X$ for each $n$. We prove that $\underline{0}\in \overline{A_\beta}$. Let $[B,\varepsilon]^s(\underline{0})$ be a neighbourhood of $\underline{0}$, where $B\in \bk$ and $\varepsilon>0$. Choose a $n_1\in \nb$  such that $\frac{1}{n_1}<\varepsilon$. Since $\uc_{\beta,n_1}$ is an open $\bs$-cover of $X$, there exist a $f\in A_{\beta}$ and a $\delta>0$ such that $B^\delta\subseteq f^{-1}(-\frac{1}{n_1},\frac{1}{n_1})$. Thus $f\in [B,\varepsilon]^s(\underline{0})\cap A_{\beta}$ and $\underline{0}\in \overline{A_\beta}=A_\beta$. Hence $\underline{0}\in A$. This completes the proof.
\ep

\section{Some observations on {\rm k}-spaces}
The concept of {\rm k}-space was studied in \cite{McCoy, McCoy-2} by R. A. McCoy in function spaces endowed with the point-open topology and the compact-open topology. For a Tychonoff space $X$ it was known that $C_p(X)$ is a {\rm k}-space if and only if $X$ satisfies $\S1(\Omega,\Gamma)$ (see \cite[Theorem 1]{McCoy}). In this section we present some observations on {\rm k}-spaces in the context of the topology of strong uniform convergence on $\bk$.
We first present the following lemma.

\begin{Lemma}
\label{Lk-1}
Let $\bk$ be a bornology with closed base on a metric space $X$. If $K\subseteq (C(X),\tau^s_\bk)$ is compact then for every $x\in X$, $\{f(x):f\in K\}$ is bounded in $\rb$.
\end{Lemma}
\bp
Fix a $x_0\in X$. Let $B=\{x_0\}$. Clearly $K\subseteq \cup\{[B,n]^s(\underline{0}):n\in \nb\}$. Since $K$ is compact, there are $n_1,\dots, n_k$ such that $K\subseteq \cup\{[B,n_i]^s(\underline{0}):i=1,\dots,k\}$. Choose $m=\max\{n_1,\dots,n_k\}$. Now for each $f\in K$, $|f(x_0)|<m$. Hence $\{f(x_0):f\in K\}$ is bounded.
\ep

Recall that with respect to the topology of strong uniform convergence on a bornology $\bk$, $C(X)$ is $\FU$ if and only if $X$ satisfies $\S1(\obs,\Gbs)$ (\cite[Corollary 2.10]{cmk}). Thus we have the following.

\begin{Th}
\label{Tk-1}
Let $\bk$ be a bornology with closed base on a metric space $X$. The following statements are equivalent.\\
$(1)$ $(C(X),\tau^s_\bk)$ is a {\rm k}-space.\\
$(2)$ $(C(X),\tau^s_\bk)$ is $\FU$.\\
$(3)$ $X$ satisfies $\S1(\obs,\Gbs)$.
\end{Th}
\bp
%The equivalence $(2)\Leftrightarrow (3)$ follows from \cite[Corollary 2.10]{cmk} and the implication $(2)\Rightarrow (1)$ is easily followed. We prove $(1)\Rightarrow (3)$.
$(1)\Rightarrow (3)$. We prove this by contrapositive argument. Suppose that $X$ does not satisfy $\S1(\obs,\Gbs)$. Then we get a sequence $\{\uc_n:n\in \nb\}$ of open $\bs$-covers of $X$ for which $\{U_n:n\in \nb\}$ is not a $\gbs$-cover of $X$ for any $U_n\in \uc_n$, $n\in \nb$. Construct a new sequence $\{\vc_n:n\in \nb\}$ of open $\bs$-covers of $X$ as follows. Define $\vc_1=\uc_1$ and for $n>1$ $\vc_n$ is a refinement of both $\vc_{n-1}$ and $\uc_n$.

For each $B\in \bk$ and each $n\in \nb$ there are a $V\in \vc_n$ and a $\delta>0$ satisfying $B^{2\delta}\subseteq V$. Let $\vc_{n,B}=\{V\in \vc_n:B^{2\delta}\subseteq V \}$. For every $V\in \vc_{n,B}$ choose a continuous function $f_{B,V}$ from $X$ to $[\frac{1}{n},n]$ such that $f_{B,V}(B^\delta)=\{\frac{1}{n}\}$ and $f_{B,V}(X\setminus V)=\{n\}$. Let $F_n=\{f_{B,V}:B\in \bk, V\in \vc_{n,B}\}$. Choose $F=\cup_{n\in \nb}F_n$ and $F^*=F\setminus \{\underline{0}\}$. It is easy to see that $\underline{0}\in \overline{F}$. Clearly $F^*$ is not closed. We now prove that $(C(X),\tau^s_\bk)$ is not a {\rm k}-space by showing that $F^*\cap K$ closed in $K$ for each compact subset $K$ of $(C(X),\tau^s_\bk)$.
Let $K$ be a compact subset of $(C(X),\tau^s_\bk)$. By Lemma \ref{Lk-1}, $\{f(x):f\in K\}$ is bounded for every $x\in X$. Define $M(x)=\sup\{f(x):f\in K\}$. Let $X_m=\{x\in X:M(x)\leq m\}$. Clearly $X_m\subseteq X_{m+1}$ for $m\in \nb$ and the collection $\{X_m:m\in \nb\}$ is a $\bs$-cover of $X$.

Suppose that for every $m,n\in \nb$, there exist a $k\geq n$ and a $V\in \vc_k$ such that $X_m\subseteq V$. We now use induction to choose a sequence $\{U_n:n\in \nb\}$, where $U_n\in \uc_n$ for each $n$. Choose a $k_1\geq 1$ and a $V_1\in \vc_{k_1}$ such that $X_1\subseteq V_1$. By construction of $\vc_{k_1}$, for each $i\in \{1,\dots,k_1\}$ choose a $U_i\in \uc_i$ such that $V_1\subseteq U_i$. Suppose that $k_m$ and $U_1,\dots, U_{k_m}$ have been chosen. Now choose a $k_{m+1}\geq k_m+1$ and a $V_{m+1}\in \vc_{k_{m+1}}$ such that $X_{m+1}\subseteq V_{m+1}$. For each $i\in \{k_m+1, \dots,k_{m+1}\}$ choose a $U_i\in \uc_i$ such that $V_{m+1}\subseteq U_i$. Proceeding in this way we obtain a sequence $\{U_n:n\in \nb\}$, where $U_n\in \uc_n$ for each $n\in \nb$, which is not a $\gbs$-cover of $X$ by our assumption. Let $B\in \bk$. There are $X_m$ and $\delta>0$ such that $B^\delta\subseteq X_m$. For $n\geq k_m$ choose a $j\geq m$ with $k_{j-1}+1\leq n\leq k_j$ and define $\delta_n=\delta$. Clearly $B^{\delta_n}\subseteq X_m\subseteq X_j\subseteq V_j\subset U_n$. Thus $B^{\delta_n}\subseteq U_n$ for all $n\geq k_m$. This shows that $\{U_n:n\in \nb\}$ is a $\gbs$-cover which is a contradiction. Hence there exist $m_0,n_0\in \nb$ such that for any $k\geq n_0$ and any $V\in \vc_k$, $X_{m_0}\nsubseteq V$.

Choose $N=\max\{m_0,n_0\}$. Consider the neighbourhood $[B,\frac{1}{N}]^s(\underline{0})$. Let $f\in [B,\frac{1}{N}]^s(\underline{0})\cap F$. Clearly $f\in F_k$ for some $k\in \nb$ and $f=f_{B,V}$, where $V\in \vc_{k,B}$. Now $\frac{1}{k}\leq f_{B,V}(x)\leq \frac{1}{N}$ for all $x\in B^\delta$, where $\delta>0$. Hence $k>N\geq n_0$ and $X_{m_0}\nsubseteq V$ for any $V\in \vc_k$. Let $x\in X_{m_0}\setminus V$. Observe that $f_{B,V}(x)=k$. Also since $k>N\geq m_0$ and $x\in X_{m_0}$, we have $k\geq M(x)$. Therefore $f_{B,V}\not\in K$ and $[B,\frac{1}{N}]^s(\underline{0})\cap F\cap K=\emptyset$. This shows that $\underline{0}$ is not a limit point of $F^*\cap K$ in $K$ and so $F^*\cap K$ must be closed in $K$. Thus $(C(X),\tau^s_\bk)$ is not a {\rm k}-space.
\ep

In combination with \cite[Theorem 3.6]{pddcsd-3},\cite[Theorem 3.2]{dcpdsd}, \cite[Proposition 2.4]{dcpdsd-2} and \cite[Corollary 2.10]{cmk}, we obtain the following.

\begin{Cor}
\label{Ckspace-1}
Let $\bk$ be a bornology with closed base on a metric space $X$. The following statements are equivalent.\\
$(1)$ $(C(X),\tau^s_\bk)$ is $\FU$.\\
$(2)$ $(C(X),\tau^s_\bk)$ is $\SFU$.\\
$(3)$ $(C(X),\tau^s_\bk)$ is $\FUf$.\\
$(4)$ $(C(X),\tau^s_\bk)$ is a {\rm k}-space.\\
$(5)$ $X$ satisfies $\S1(\obs,\Gbs)$.\\
$(6)$ $X$ satisfies $\Sf(\obs,\Gbs)$.\\
$(7)$ $X$ satisfies $\binom{\obs}{\Gbs}$.\\
$(8)$ $X$ satisfies $\Split(\obs,\Gbs)$.
\end{Cor}

\begin{Cor}[{\rm see \cite[Theorem 1]{McCoy-2}}]
\label{Ckk-1}
For a Tychonoff space $X$ the following statements are equivalent.\\
$(1)$ $C_k(X)$ is a {\rm k}-space.\\
$(2)$ $C_k(X)$ is $\FU$.\\
$(3)$ $X$ satisfies $\S1(\kc,\Gamma_k)$.
\end{Cor}

As an immediate application of the preceding result, we present two examples one of which is a {\rm k}-space and other is not a {\rm k}-space.
\begin{Ex}
\label{Ek-1}
{\rm Let $X=\rb$ and $\bk$ be a bornology generated by $\{(-x,x):x>0\}$. Let $\{\uc_n:n\in \nb\}$ be a sequence of open $\bs$-covers of $X$. For each $n\in \nb$, choose a $U_n\in \uc_n$ such that $(-n,n)\subseteq U_n$. It is easy to see that $\{U_n:n\in \nb\}$ is a $\gbs$-cover of $X$. By Theorem \ref{Tk-1}, $(C(X),\tau^s_\bk)$ is a {\rm k}-space.}
\end{Ex}

\begin{Ex}
\label{Ek-2}
{\rm Consider the Baire space $X=\nb^\nb$ endowed with the Baire metric $\rho$ and the bornology $\bk=\fc$. We first show that $X$ does not satisfy $\S1(\obs,\Gbs)$. To see this, consider $\uc_n=\{U^n_m:m\in \nb\}$, where $U^n_m=\{f\in X:f(n)\leq m\}$. Let there exist a $U^n_{m_n}\in \uc_n$ for each $n\in \nb$ such that $\{U^n_{m_n}:n\in \nb\}$ is a $\gbs$-cover of $X$. Choose an $h\in X$ such that $h(n)=m_n+1$ for each $n\in \nb$ and a $B\in \bk$ with $h\in B$. Since $\{U^n_{m_n}:n\in \nb\}$ is a $\gbs$-cover of $X$, for $B\in \bk$ there exist a $n_0\in \nb$ and a sequence $\{\delta_n:n\geq n_0\}$ of positive real numbers satisfying $B^{\delta_n}\subseteq U^n_{m_n}$ for all $n\geq n_0$. It follows that $h(n)\leq m_n$ for all $n\geq n_0$, which is a contradiction. Hence $X$ does not satisfy $\S1(\obs,\Gbs)$. By Theorem \ref{Tk-1}, $(C(X),\tau^s_\bk)$ is not a {\rm k}-space.}
\end{Ex}

\begin{Prop}
\label{Pk}
Let $\bk$ be a bornology with closed base on a metric space $X$. If $\bk$ has a countable base, then $(C(X),\tau^s_\bk)$ is a {\rm k}-space.
\end{Prop}
\bp
Let $\bk_0=\{B_n:n\in \nb\}$ be a base for $\bk$. Without loss of generality we assume that $B_n\subseteq B_{n+1}$ for each $n$. Let $A\subseteq C(X)$ be such that $A\cap K$ is closed in $K$ for every compact subset $K$ of $C(X)$. We show that $A$ is closed. Let $\underline{0}\in \overline{A}$. Choose a $f_n\in A\cap [B_n,\frac{1}{n}]^s(\underline{0})$ for each $n\in \nb$. We claim that the sequence $\{f_n:n\in \nb\}$ converges to $\underline{0}$. Let $[B,\varepsilon]^s(\underline{0})$ be a neighbourhood of $\underline{0}$, where $B\in \bk$ and $\varepsilon>0$. Choose a $n_0\in \nb$ such that $B\subseteq B_n$ and $\frac{1}{n}<\varepsilon$ for all $n\geq n_0$. Clearly $f_n\in [B,\varepsilon]^s(\underline{0})$ for all $n\geq n_0$. Therefore $\{f_n:n\in \nb\}$ converges to $\underline{0}$. Let $K=\{\underline{0}\}\cup \{f_n:n\in \nb\}$. Clearly $\underline{0}\in \overline{A\cap K}$. Since $A\cap K$ is closed, $\underline{0}\in A$. Hence $A$ is closed.
\ep

\begin{Prop}
\label{Pk-A}
Let $\bk$ be a bornology with closed base $\bk_0$ on a metric space $X$ and $X$ be $\bs$-Lindel\"{o}f. If $|\bk_0|<\pf$, then $(C(X),\tau^s_\bk)$ is a {\rm k}-space.
\end{Prop}
\bp
By \cite[Theorem 3.1]{dcpdsd-2}, $X$ is a $\gbs$-set. Therefore $X$ satisfies $\S1(\obs,\Gbs)$. By Theorem \ref{Tk-1}, $(C(X),\tau^s_\bk)$ is a {\rm k}-space.
\ep

\begin{Prop}
\label{Pk-1}
Let $\bk$ be a bornology with a closed base $\bk_0$ which contains the nonempty compact subsets of a metric space $X$. If $(C(X),\tau^s_\bk)$ is a {\rm k}-space, then $\bk$ is local.
\end{Prop}
\bp
Assume that $\bk$ is not local. Choose a $x_0\in X$ such that $S_\delta(x_0)\setminus B\neq \emptyset$ for any $B\in \bk$ and any $\delta>0$. Consider a sequence $\{S_\frac{1}{n}(x_0):n\in \nb\}$ of neighbourhoods of $x_0$. Let $x_{n,B}\in S_\frac{1}{n}(x_0)\setminus B$ for each $B\in \bk_0$ and each $n\in \nb$. Since $B$ is closed, $d(x_{n,B},B)>0$. Choose $\delta_{n,B}=\frac{1}{2}d(x_{n,B},B)$. Clearly $S_\frac{1}{n}(x_0)\setminus B^{\delta_{n,B}}\neq \emptyset$. For each $n\in \nb$, consider the collection $\uc_n=\{B^{\delta_{n,B}}:B\in \bk_0\}$, which is an open $\bs$-cover of $X$. By the given condition, there exists a $U_n\in \uc_n$ for each $n$ such that $\{U_n:n\in \nb\}$ is a $\gbs$-cover of $X$. Let $y_n\in S_\frac{1}{n}(x_0)\setminus U_n$. Clearly the sequence $\{y_n:n\in \nb\}$ converges to $x_0$. Now $A=\{x_0\}\cup \{y_n:n\in \nb\}$ is compact and so $A\in \bk_0$. Clearly $A\nsubseteq U_n$ for any $n\in \nb$. This a contradiction to the fact that $\{U_n:n\in \nb\}$ is a $\gbs$-cover of $X$. Hence $\bk$ is local.
\ep

\begin{Prop}
\label{Pk-2}
Let $\bk$ be a bornology with closed base on metric space $X$. If $(C(X),\tau^s_\bk)$ is a {\rm k}-space, then it has countable tightness.
\end{Prop}
\bp
By the given condition $X$ satisfies $\S1(\obs,\Gbs)$. Also by \cite[Theorem 2.8]{cmk}, we obtain that every open $\bs$-cover of $X$ has a $\gbs$-cover of $X$. Hence $X$ is $\bs$-Lindel\"{o}f. Now by \cite[Theorem 2.1]{cmk}, $(C(X),\tau^s_\bk)$ has countable tightness.
\ep

Using Theorem \ref{Tk-1}, \cite[Lemma 5.3]{dcpdsd-2} and \cite[Proposition 2.1]{dcpdsd-2} we have the following.
\begin{Prop}
\label{Pk-3}
Let $\bk$ be a bornology with compact base on a metric space $X$. \\
$(1)$ Let $(C(X),\tau^s_\bk)$ be a {\rm k}-space and $Y\subseteq X$ be closed. Then $(C(Y),\tau^s_{\bk_Y})$ is a {\rm k}-space, where $\bk_Y=\{B\cap Y:B\in \bk\}$.\\
$(2)$ Let $f:X\rightarrow Y$ be a continuous function. If $(C(X),\tau^s_\bk)$ is a {\rm k}-space, then $(C(Y),\tau^s_{f(\bk)})$ is a {\rm k}-space.
\end{Prop}

\section{Topological games and Discrete Selectivity}

\subsection{Certain observations on strong $\bk$-open game}

We first introduce two games.
\begin{Def}
\label{Dobsx}
{\rm The} strong $\bk$-open {\rm (in short,} $\bs$-open {\rm)  game on $X$ is defined as follows. Suppose that an infinitely long game is played by two players ONE and TWO. In the $n$th round, ONE chooses $B_n$, where $B_n\in \bk$, and TWO responds by choosing an open set $U_n$ satisfying $B_n^{\delta_n}\subseteq U_n$ for some $\delta_n>0$. ONE wins the play if $\{U_n:n\in \nb\}$ is an open $\bs$-cover of $X$. Otherwise TWO wins.}
\end{Def}
We denote the $\bs$-open game on $X$ by $\obs(X)$.

\begin{Def}
\label{Dgbsx}
{\rm The} $\gbs$-open {\rm game on $X$ is defined as follows. Suppose that an infinitely long game is played by two players ONE and TWO. In the $n$th round, ONE chooses $B_n$, where $B_n\in \bk$, and TWO responds by choosing an open set $U_n$ satisfying $B_n^{\delta_n}\subseteq U_n$ for some $\delta_n>0$. ONE wins the play if $\{U_n:n\in \nb\}$ is a $\gbs$-cover of $X$. Otherwise TWO wins.}
\end{Def}
We denote the $\gbs$-open game on $X$ by $\Gbs(X)$.

To prove the following result we use the technique of \cite[Lemma 4.1]{mn} and the fact that a sequence of open sets $\{U_n:n\in \nb\}$ is a $\gbs$-cover of $X$ if and only if every infinite subsequence of $\{U_n:n\in \nb\}$ is an open $\bs$-cover of $X$.

\begin{Th}
\label{Ttg-1}
Let $\bk$ be a bornology with closed base on a metric space $X$. The following statements are equivalent.\\
%$(1)$ ONE has a winning strategy in the game $\obs(X)$.\\
%$(2)$ ONE has a winning strategy in the game $\Gbs(X)$.
$(1)$ ONE has a $\WS$ in $\obs(X)$.\\
$(2)$ ONE has a $\WS$ in $\Gbs(X)$.
\end{Th}
\bp
$(1)\Rightarrow (2)$. Let $\psi$ be a $\WS$ for ONE in $\obs(X)$. We define a strategy $\sigma$ for ONE in $\Gbs(X)$ as follows. Define $\sigma(\emptyset)=\psi(\emptyset)$, where $\psi(\emptyset)\in \bk$. Let TWO respond by choosing an open set $U_1$ such that $\sigma(\emptyset)^{\delta_1}\subseteq U_1$ for some $\delta_1>0$. Suppose that $U_1, \dots, U_{n-1}$ have already been chosen. We now define $\sigma(U_1,\dots,U_{n-1})$ in such a way that  $\psi(\emptyset)\cup [\cup\{\psi(U_{i_1},\dots U_{i_k}): 1\leq i_1\leq \cdots \leq i_k\leq n-1\}]\subseteq \sigma(U_1,\dots,U_{n-1})$. Let TWO respond by choosing an open set $U_n$ such that $\sigma(U_1,\dots,U_{n-1})^{\delta_n}\subseteq U_n$ for some $\delta_n>0$. This defines the strategy $\sigma$ for ONE in $\Gbs(X)$. We show that the sequence $\{U_n:n\in \nb\}$ of moves of TWO in $\Gbs(X)$ forms a $\gbs$-cover of $X$. For this we prove that every infinite subsequence of $\{U_n:n\in \nb\}$ is an open $\bs$-cover of $X$. Let $\{U_{n_k}:k\in \nb\}$ be an infinite subsequence of $\{U_n:n\in \nb\}$. Clearly $\psi(\emptyset)\subseteq \sigma(U_1,\dots, U_{{n_1}-1})$. Since $\sigma(U_1,\dots, U_{{n_1}-1})^{\delta_{n_1}}\subseteq U_{n_1}$ for some $\delta_{n_1}>0$, $\psi(\emptyset)^{\delta_{n_1}}\subseteq U_{n_1}$. Again $\psi(U_{n_1},\dots, U_{n_k})\subseteq \sigma(U_1,\dots, U_{n_{k+1}-1})$ and since $\sigma(U_1,\dots, U_{n_{k+1}-1})^{\delta_{n_{k+1}}}\subseteq U_{n_{k+1}}$ for some $\delta_{n_{k+1}}>0$, $\psi(U_{n_1},\dots, U_{n_k})^{\delta_{n_{k+1}}}\subseteq U_{n_{k+1}}$ and so on.

Now $\psi(\emptyset), U_{n_1},\dots, \psi(U_{n_1},\dots,U_{n_k}), U_{n_{k+1}},\dots$ is a legitimate play in $\obs(X)$. Now $\{U_{n_k}:k\in \nb\}$ is an open $\bs$-cover of $X$ as $\psi$ is a $\WS$ for ONE in $\obs(X)$. Therefore $\{U_n:n\in \nb\}$ is a $\gbs$-cover of $X$. Hence $\sigma$ is a $\WS$ for ONE in $\Gbs(X)$.
\ep

For the next two results we use the method of Theorems $15,17$ given in \cite{ch}

\begin{Th}
\label{Ttg-2}
Let $\bk$ be a bornology with closed base on a metric space $X$. The following statements are equivalent.\\
$(1)$ ONE has a $\WS$ in $\obs(X)$.\\
$(2)$ ONE has a $\WS$ in $\Gbs(X)$.\\
$(3)$ TWO has a $\WS$ in $\G1(\obs,\obs)$ on $X$.\\
$(4)$ TWO has a $\WS$ in $\G1(\obs,\Gbs)$ on $X$.
\end{Th}
\bp
We prove $(1)\Leftrightarrow (3)$.

$(1)\Rightarrow (3)$. Let $\psi$ be a $\WS$ for ONE in $\obs(X)$. Define a strategy $\sigma$ for TWO in $\G1(\obs,\obs)$ as follows. Let the first move of ONE in $\G1(\obs,\obs)$ be $\uc_1$. For $\psi(\emptyset)\in \bk$ there are a $U_1\in \uc_1$ and a $\delta_1>0$ satisfying $\psi(\emptyset)^{\delta_1}\subseteq U_1$ for some $\delta_1>0$. Define $\sigma(\uc_1)=U_1$. Suppose that $U_1,\dots,U_{n-1}$ have been chosen. In the $n$th round, let the move of ONE in $\G1(\obs,\obs)$ be $\uc_n$. For $\psi(U_1,\dots,U_{n-1})\in \bk$ there are a $U_n\in \uc_n$ and a $\delta_n>0$ satisfying $\psi(U_1,\dots U_{n-1})^{\delta_n}\subseteq U_n$. Define $\sigma(\uc_1,\dots,\uc_n)=U_n$.

For a play $\uc_1,\sigma(\uc_1),\dots,\uc_n,\sigma(\uc_1,\dots,\uc_n),\dots$ in $\G1(\obs,\obs)$, the corresponding play in $\obs(X)$ is $\psi(\emptyset), U_1, \dots, \psi(U_1,\dots, U_{n-1}), U_n, \dots$. Since $\psi$ is a $\WS$ for ONE in $\obs(X)$, $\{U_n:n\in \nb\}$ is an open $\bs$-cover of $X$. Hence $\sigma$ is a $\WS$ for TWO in $\G1(\obs,\obs)$.

$(3)\Rightarrow (1)$. Let $\sigma$ be a $\WS$ for TWO in $\G1(\obs,\obs)$. Define a strategy $\psi$ for ONE in $\obs(X)$ as follows. Let $\{U_1,\dots,U_{n-1}\}$ be a finite sequence of open subsets of $X$. Assume that $\psi(U_1,\dots,U_k)$ has been defined for all $k<n-1$ and $\uc_k$, an open $\bs$-cover of $X$, is defined for all $k\leq n-1$.

Suppose that for each $B\in \bk$ there exist an open set $U_B$ and a $\delta>0$ satisfying $B^\delta\subseteq U_B$ such that for every open $\bs$-cover $\uc$, $\sigma(\uc_1,\dots,\uc_{n-1},\uc)\neq U_B$. Clearly $\uc'=\{U_B:B\in \bk\}$ is an open $\bs$-cover of $X$. Now $\sigma(\uc_1,\dots,\uc_{n-1},\uc')=U_B$ for some $B\in \bk$, which is a contradiction. Hence there exists a $B_n\in \bk$ such that for every open set $U$ with $B^\delta_n\subseteq U$ for some $\delta>0$ there exists an open $\bs$-cover $\uc$ satisfying $\sigma(\uc_1,\dots, \uc_{n-1},\uc)=U$. We define $\psi(U_1,\dots, U_{n-1})=B_n$. Thus $\psi$ is a strategy for ONE in $\obs(X)$.

In response to the move $\psi(U_1,\dots, U_{n-1})$ of ONE in $\obs(X)$, let TWO choose $U_n$, where $B^{\delta_n}_n\subseteq U_n$ for some $\delta_n>0$. Now for $U_n$ there
exists an open $\bs$-cover $\uc_n$ of $X$ such that $\sigma(\uc_1,\dots, \uc_{n-1},\uc_n)=U_n$. For a play in $\obs(X)$ $$\psi(\emptyset), U_1,\dots,\psi(U_1,\dots,U_{n-1}), U_n,\dots,$$
the corresponding play in $\G1(\obs,\obs)$ is $$\uc_1,\sigma(\uc_1),\dots,\uc_n,\sigma(\uc_1,\dots,\uc_n),\dots.$$
Since $\sigma$ is a $\WS$ for TWO in $\G1(\obs,\obs)$, $\{U_n:n\in \nb\}$ is an open $\bs$-cover of $X$. Hence $\psi$ is a $\WS$ for ONE in $\obs(X)$.

\ep

%To prove the following result we use the technique of \cite[Theorem 17]{ch}.
\begin{Th}
\label{Ttg-3}
Let $\bk$ be a bornology with closed base on a metric space $X$. The following statements are equivalent.\\
$(1)$ ONE has a $\PWS$ in $\obs(X)$.\\
$(2)$ ONE has a $\PWS$ in $\Gbs(X)$.\\
$(3)$ TWO has a $\WMS$ in $\G1(\obs,\obs)$ on $X$.\\
$(4)$ TWO has a $\WMS$ in $\G1(\obs,\Gbs)$ on $X$.
\end{Th}
\bp
%We prove only $(1)\Leftrightarrow (3)$.

$(1)\Rightarrow (3)$. Let $\psi$ be a $\PWS$ for ONE in $\obs(X)$. We define a $\MS$ $\sigma$ for TWO in $\G1(\obs,\obs)$ as follows. Let $n\in \nb$ and $\uc$ be an open $\bs$-cover of $X$. For $\psi(n)$ there exist a $\delta>0$ and a $U\in \uc$ such that $\psi(n)^\delta\subseteq U$. Define $\sigma(\uc,n)=U$. Thus $\sigma$ is a $\MS$ for TWO in $\G1(\obs,\obs)$. For a play $\uc_1, \sigma(\uc_1,1),\dots,\uc_n,\sigma(\uc_n,n)\dots$ in $\G1(\obs,\obs)$, the corresponding play in $\obs(X)$ is $\psi(1), U_1,\dots,\psi(n),U_n\dots$, where $\sigma(\uc_n,n)=U_n$ for $n\in \nb$. Since $\psi$ is a $\WS$ for ONE, $\{U_n:n\in \nb\}$ is an open $\bs$-cover of $X$. Hence $\sigma$ is a $\WS$ for TWO in $\G1(\obs,\obs)$.

$(3)\Rightarrow (1)$. Let $\sigma$ be a $\WMS$ for TWO in $\G1(\obs,\obs)$. We define a $\PS$ $\psi$ for ONE in $\obs(X)$ as follows. Fix a $n\in \nb$. Suppose that for each $B\in \bk$ there exists an open set $U_B$ with $B^\delta\subseteq U_B$ such that for every open $\bs$-cover $\uc$, $\sigma(\uc,n)\neq U_B$. Now $\sigma(\{U_B:B\in \bk\},n)\not\in \{U_B:B\in \bk\}$, a contradiction. Hence there exists a $B_n$ such that for every open set $U$ with $B_n^\delta\subseteq U$, there exists an open $\bs$-cover $\uc_n$ satisfying $\sigma(\uc_n,n)=U$. Define $\psi(n)=B_n$. Thus $\psi$ is a $\PS$ for ONE in $\obs(X)$. For $n\in \nb$ and $\psi(n)$, let TWO's response in $\obs(X)$ is $U_n$. Now choose an open $\bs$-cover $\uc_n$ satisfying $\sigma(\uc_n,n)=U_n$. For a play $\psi(1), U_1,\dots,\psi(n),U_n,\dots$ in $\obs(X)$, the corresponding play in $\G1(\obs,\obs)$ is $\uc_1,\sigma(\uc_1,1),\dots,\uc_n,\sigma(\uc_n,n),\dots$. Since $\sigma$ is a $\WS$ for TWO, $\{U_n:n\in \nb\}$ is an open $\bs$-cover of $X$. Hence $\psi$ is a $\WS$ for ONE in $\obs(X)$.

Other implications can be similarly verified.
\ep

We skip the proof of the next result.
\begin{Prop}
\label{Ptg-1}
Let $\bk$ be a bornology with closed base on a metric space $X$. If TWO has a $\WS$ in $\obs(X)$, then TWO has a $\WS$ in $\Gbs(X)$.
\end{Prop}

\begin{Th}
\label{TTT}
Let $\bk$ be a bornology with closed base on a metric space $X$. The following statements are equivalent.\\
%$(1)$ TWO has a winning strategy in the game $\obs(X)$.\\
%$(2)$ ONE has a winning strategy in the game $\G1(\obs,\obs)$ on $X$.
$(1)$ TWO has a $\WS$ in $\obs(X)$.\\
$(2)$ ONE has a $\WS$ in $\G1(\obs,\obs)$ on $X$.
\end{Th}
\bp
$(1)\Rightarrow (2)$. Let $\psi$ be a $\WS$ for TWO in $\obs(X)$. We define a strategy $\sigma$ for ONE in $\G1(\obs,\obs)$ as follows. For a finite sequence $t$ of members from $\bk$ let $\uc_t=\{\psi(t\frown<B>):B\in \bk\}$. Define $\sigma(\emptyset)=\uc_\emptyset$. Let TWO respond with $U_1\in \sigma(\emptyset)$, where $U_1=\psi(B_1)$. Suppose that $U_1, \dots, U_{n-1}$ have been chosen. Now define $\sigma(U_1, \dots, U_{n-1})=\uc_{<B_1,\dots,B_{n-1}>}$. Let TWO respond with $U_n$, where $U_n=\psi(B_1,\dots,B_n)$. A play in $\G1(\obs,\obs)$ is $$\sigma(\emptyset),U_1,\dots,\sigma(U_1,\dots,U_{n-1}),U_n,\dots.$$
The corresponding play in $\obs(X)$ is $$B_1,\psi(B_1),\dots,B_n,\psi(B_1,\dots,B_n)\dots.$$
Since $\psi$ is a $\WS$ for TWO in $\obs(X)$, $\{U_n:n\in \nb\}$ is not an open $\bs$-cover of $X$. Hence $\sigma$ is a $\WS$ for ONE in $\G1(\obs,\obs)$.

$(2)\Rightarrow (1)$. Let $\sigma$ be a $\WS$ for ONE in $\G1(\obs,\obs)$. Define a strategy $\psi$ for TWO in $\obs(X)$ as follows. Let $B_1$ be the first move of ONE in $\obs(X)$. For $B_1\in \bk$ there exist a $U_1\in \sigma(\emptyset)$ and a $\delta_1>0$ such that $B_1^{\delta_1}\subseteq U_1$. Define $\psi(B_1)=U_1$. In $\G1(\obs,\obs)$, let TWO respond with $U_1$. Suppose that $U_1,\dots, U_{n-1}$ have been chosen. Now define $\psi(B_1,\dots,B_n)=U_n$, where $B_n^{\delta_n}\subseteq U_n$ for some $\delta>0$ and $U_n\in \sigma(U_1,\dots,U_{n-1})$. In $\G1(\obs,\obs)$, let TWO respond with $U_n$ and so on. A play in $\obs(X)$ is $$B_1,\psi(B_1), \dots,B_n,\psi(B_1,\dots, B_n),\dots$$ The corresponding play in $\G1(\obs,\obs)$ is $$\sigma(\emptyset),U_1 \dots,\sigma(U_1,\dots,U_{n-1}),U_n,\dots$$
Since $\sigma$ is a $\WS$ for ONE in $\G1(\obs,\obs)$, $\{U_n:n\in \nb\}$ is not an open $\bs$-cover of $X$. Hence $\psi$ is a $\WS$ for TWO in $\obs(X)$.
\ep

%The following result is obtained similarly.
\begin{Th}
\label{TTT1}
Let $\bk$ be a bornology with closed base on a metric space $X$. The following statements are equivalent.\\
%$(1)$ TWO has a winning strategy in the game $\Gbs(X)$.\\
%$(2)$ ONE has a winning strategy in the game $\G1(\obs,\Gbs)$ on $X$.
$(1)$ TWO has a $\WS$ in $\Gbs(X)$.\\
$(2)$ ONE has a $\WS$ in $\G1(\obs,\Gbs)$ on $X$.
\end{Th}

\begin{Th}
\label{TT}
Let $\bk$ be a bornology with closed base on a metric space $X$. The following statements are equivalent.\\
%$(1)$ TWO has a winning Markov strategy in the game $\obs(X)$.\\
%$(2)$ ONE has a predetermined winning strategy in the game $\G1(\obs,\obs)$ on $X$.
$(1)$ TWO has a $\WMS$ in $\obs(X)$.\\
$(2)$ ONE has a $\PWS$ in $\G1(\obs,\obs)$ on $X$.
\end{Th}
\bp
$(1)\Rightarrow (2)$. Let $\psi$ be a $\WMS$ for TWO in $\obs(X)$. We define a $\PS$ $\sigma$ for ONE in $\G1(\obs,\obs)$ as follows. Fix a $n\in \nb$. For $B\in \bk$ $\psi(B,n)$ is an open set satisfying $B^\delta\subseteq \psi(B,n)$ for some $\delta>0$. Thus $\{\psi(B,n):B\in \bk\}$ is an open $\bs$-cover of $X$. Define $\sigma(n)=\{\psi(B,n):B\in \bk\}$. Hence $\sigma$ is a $\PS$ for ONE in $\G1(\obs,\obs)$. It is easy to see that $\sigma$ is a $\WS$.
%For a play $\sigma(1), U_1,\dots,\sigma(n), U_n,\dots$ in $\G1(\obs,\obs)$, where $U_n=\psi(B_n,n)$, $B_n\in \bk$, the corresponding play in $\obs(X)$ is $B_1,\psi(B_1,1),\dots,B_n,\psi(B_n,n),\dots$. Since $\psi$ is a $\WS$ for TWO in $\obs(X)$, $\{U_n:n\in \nb\}$ is not an open $\bs$-cover of $X$. Hence $\sigma$ is a $\WS$ for ONE in $\G1(\obs,\obs)$.

$(2)\Rightarrow (1)$. Let $\sigma$ be a $\PWS$ for ONE in $\G1(\obs,\obs)$. We define a $\MS$ $\psi$ for TWO in $\obs(X)$ as follows. For $B\in \bk$ and $n\in \nb$ there exists a $U\in \sigma(n)$ and a $\delta>0$ such that $B^\delta\subseteq U$. Define $\psi(B,n)=U$. Thus $\psi$ is a $\MS$ for TWO in $\obs(X)$. It is easy to see that $\psi$ is a $\WS$.
%For a play $B_1,\psi(B_1,1),\dots,B_n,\psi(B_n,n),\dots$ in $\obs(X)$, the corresponding play in $\G1(\obs,\obs)$ is $\sigma(1), U_1,\dots,\sigma(n),U_n,\dots$, where $U_n=\psi(B_n,n)$ for $n\in \nb$. Since $\sigma$ is a $\WS$ for ONE in $\G1(\obs,\obs)$, $\{U_n:n\in \nb\}$ is not an open $\bs$-cover of $X$. Hence $\psi$ is a $\WS$ for TWO in $\obs(X)$.
\ep

%The following result is obtained similarly.
\begin{Th}
\label{TT1}
Let $\bk$ be a bornology with closed base on a metric space $X$. The following statements are equivalent.\\
%$(1)$ TWO has a winning Markov strategy in $\Gbs(X)$.\\
%$(2)$ ONE has a predetermine winning strategy in $\G1(\obs,\Gbs)$ on $X$.
$(1)$ TWO has a $\WMS$ in $\Gbs(X)$.\\
$(2)$ ONE has a $\PWS$ in $\G1(\obs,\Gbs)$ on $X$.
\end{Th}

\subsection{Results on discrete selectivity and related games}
The concept of discrete selectivity and related games were introduced and studied in \cite{tkachuk-1,tkachuk-2} by V.V. Tkachuk in function spaces with respect to the point-open topology. $C_p(X)$ is discretely selective if and only if $X$ is uncountable, provided $X$ is a Tychonoff space \cite[Proposition 3.3]{tkachuk-1}. In this section we study discretely selective property and related games in function spaces with respect to the topology of strong uniform convergence on a bornology.

\begin{Th}
\label{Tds-1}
Let $\bk$ be a bornology with a compact base $\bk_0$ on a metric space $X$. The following statements are equivalent.\\
$(1)$ $\bk_0$ is uncountable.\\
$(2)$ $(C(X),\tau^s_\bk)$ is discretely selective.
\end{Th}
\bp
$(1)\Rightarrow (2)$. Let $\{U_n:n\in \nb\}$ be a sequence of nonempty open subsets of $(C(X),\tau^s_\bk)$. Without loss of generality we assume that $U_n=[B_n,\frac{1}{n}]^s(g_n)$ for some $B_n\in \bk_0$ and $g_n\in C(X)$. Since $\bk_0$ is uncountable, there exists a $B\in \bk_0\setminus \{B_n:n\in \nb\}$. For each $n$ let $x_n\in B\setminus B_n$. Since $B_n$ is closed, there exists a $\delta_n>0$ such that $S_{2\delta_n}(x_n)\cap B_n=\emptyset$. Choose a $f_n\in C(X)$ satisfying $f_n(x)=g_n(x)$ for all $x\in \overline{B^{\delta_n}}$ and $f_n(x_n)=n$. Clearly $f_n\in [B_n,\frac{1}{n}]^s(g_n)$. We now show that $\{f_n:n\in \nb\}$ $(=D)$ is closed and discrete. Let $f\in C(X)$. As $B$ is compact, $f(B)$ is bounded and let $m=\sup\{2|f(x)|:x\in B\}$. Consider the neighbourhood $[B,\frac{m}{2}]^s(f)$. If $f_n\in [B,\frac{m}{2}]^s(f)$, then $|f_n(x)|< \frac{m}{2}+|f(x)|$ for all $x\in B^\delta$, where $\delta>0$. Thus $|f_n(x_n)|< m$ and so $n< m$. Therefore $[B,\frac{m}{2}]^s(f)\cap D\subseteq \{f_1,\dots,f_{m-1}\}$. Hence $D$ is closed and discrete.

$(2)\Rightarrow (1)$. If $\bk_0$ is countable, then $(C(X),\tau^s_\bk)$ is first countable by \cite[Theorem 3.1]{cmh}. Again by $(2)$ and \cite[Proposition 3.2(b)]{tkachuk-1}, $(C(X),\tau^s_\bk)$ is discrete, which is a contradiction. Hence $\bk_0$ is uncountable.
\ep

For a Tychonoff space $X$ it is known that the player ONE has a $\WS$ in $\g(C_p(X),u)$, $\CL(C_p(X),u)$ for some $u\in C_p(X)$ and $\CD(C_p(X))$ respectively are equivalent to the fact that the player ONE has a $\WS$ in the point-open game on $X$ \cite[Theorem 3.8]{tkachuk-2}. A similar result holds for the space $(C(X),\tau^s_\bk)$.

\begin{Th}
\label{Tg-5}
Let $\bk$ be a bornology with a compact base $\bk_0$ on a metric space $X$. The following statements are equivalent.\\
%$(1)$ ONE has a winning strategy in the game $\g(\tau^s_\bk)$.\\
%$(2)$ ONE has a winning strategy in the game $\CL(\tau^s_\bk)$. \\
%$(3)$ ONE has a winning strategy in the game $\CD(\tau^s_\bk)$.\\
%$(4)$ ONE has a winning strategy in the game $\obs(X)$.
$(1)$ ONE has a $\WS$ in $\g(\tau^s_\bk)$.\\
$(2)$ ONE has a $\WS$ in $\CL(\tau^s_\bk)$. \\
$(3)$ ONE has a $\WS$ in $\CD(\tau^s_\bk)$.\\
$(4)$ ONE has a $\WS$ in $\obs(X)$.
\end{Th}
\bp
We prove only the implications $(2)\Rightarrow (3)\Rightarrow (4)\Rightarrow (1)$.

$(2)\Rightarrow (3)$ Let $\sigma$ be a $\WS$ for ONE in $\CL(\tau^s_\bk)$. We define another strategy $\psi$ for ONE in $\CL(\tau^s_\bk)$ as follows. Define $\psi(\emptyset)=\sigma(\emptyset)\setminus \{\underline{0}\}$. Suppose that $f_1, \dots, f_{n-1}$ are already chosen. Now define $\psi(f_1,\dots, f_{n-1})=\sigma(f_1,\dots f_{n-1})\setminus \{\underline{0}\}$. TWO responds by choosing a $f_n\in \psi(f_1,\dots, f_{n-1})$. Since $\sigma$ is a $\WS$ for ONE, $\underline{0}\in \overline{\{f_n:n\in \nb\}}$. Therefore $\psi$ is also a $\WS$ for ONE in $\CL(\tau^s_\bk)$. Now $\underline{0}\in \overline{\{f_n:n\in \nb\}}\setminus \{f_n:n\in \nb\}$. This shows that $\{f_n:n\in \nb\}$ is not closed and hence $\psi$ is a $\WS$ for ONE in $\CD(\tau^s_\bk)$.

$(3)\Rightarrow (4)$. Let $\psi$ be a $\WS$ for ONE in $\CD(\tau^s_\bk)$. We define a strategy $\tau$ for ONE in $\obs(X)$ as follows. Let $[B_1,\varepsilon_1]^s(g_1)\subseteq \psi(\emptyset)$ for some $B_1\in \bk$, $\varepsilon_1>0$ and $g_1\in \psi(\emptyset)$. Define $\tau(\emptyset)=B_1$. TWO responds by choosing an open set $U_1$ such that $B_1^{2\delta_1}\subseteq U_1$ for some $\delta_1>0$. Choose a $f_1\in C(X)$ satisfying $f_1(x)=g_1(x)$ for all $x\in \overline{B_1^{\delta_1}}$ and $f_1(x)=1$ for all $x\in X\setminus U_1$. Clearly $f_1\in \psi(\emptyset)$. In $\CD(\tau^s_\bk)$, let TWO respond by choosing $f_1$. Suppose that $f_1, \dots, f_{n-1}$ and $U_1, \dots, U_{n-1}$ are already chosen. Let $[B_n,\varepsilon_n]^s(g_n)\subseteq \psi(f_1,\dots,f_{n-1})$ for some $B_n\in \bk$, $\varepsilon_n>0$ and $g_n\in \psi(f_1,\dots,f_{n-1})$. Define $\tau(U_1,\dots, U_{n-1})=B_n$. TWO responds by choosing an open set $U_n$ such that $B^{2\delta_n}_n\subseteq U_n$ for some $\delta_n>0$. Choose a $f_n\in C(X)$ satisfying $f_n(x)=g_n(x)$ for all $x\in \overline{B_n^{\delta_n}}$ and $f_n(x)=n$ for all $x\in X\setminus U_n$. In $\CD(\tau^s_\bk)$, let TWO respond by choosing $f_n$ and so on. The play in $\CD(\tau^s_\bk)$ is

$\psi(\emptyset), f_1,\dots,\psi(f_1,\dots, f_{n-1}), f_n, \dots$

The corresponding play in $\obs(X)$ is

$\tau(\emptyset),U_1, \dots, \tau(U_1,\dots, U_{n-1}), U_n, \dots$

We now prove that $\{U_n:n\in \nb\}$ is an open $\bs$-cover of $X$. Assume the contrary. There exists a $B\in \bk_0$ such that $B^\delta\nsubseteq U_n$ for any $\delta>0$ and any $n\in \nb$. It follows that $B\nsubseteq U_n$ for any $n$, as $B$ is compact. Let $f\in C(X)$ and $m=\sup\{2|f(x)|:x\in B\}$. Consider the neighbourhood $[B,\frac{m}{2}]^s(f)$. It can be easily seen that $[B,\frac{m}{2}]^s(f)\cap \{f_n:n\in \nb\}\subseteq \{f_1,\dots,f_{m-1}\}$ and hence the collection $\{f_n:n\in \nb\}$ is discrete and closed. This contradicts the fact that $\psi$ is a $\WS$ for ONE in $\CD(\tau^s_\bk)$. Therefore our assumption is wrong. Hence $\tau$ is a $\WS$ for ONE in $\obs(X)$.

$(4)\Rightarrow (1)$. By Theorem \ref{Ttg-1}, let $\sigma$ be a $\WS$ for ONE in $\Gbs(X)$. We define a strategy $\psi$ for ONE in $\g(\tau^s_\bk)$ as follows.
%Define $\psi(\emptyset)=[B_1,1]^s(\underline{0})$ where $\sigma(\emptyset)=B_1\in \bk$. Let TWO choose $f_1\in \psi(\emptyset)$. There exists a $\delta_1>0$ such that $B_1^{\delta_1}\subseteq f^{-1}_1(-1,1)$. In $\Gbs(X)$, let TWO choose $f^{-1}_1(-1,1)$ ($= U_1$).
Define $\psi(f_1,\dots,f_{n-1})=[B_n,\frac{1}{n}]^s(\underline{0})$, where $\sigma(U_1,\dots, U_{n-1})=B_n$. TWO responds by choosing a $f_n\in [B_n,\frac{1}{n}]^s(\underline{0})$. There is a $\delta_n>0$ satisfying $B_n^{\delta_n}\subseteq f^{-1}_n(-\frac{1}{n},\frac{1}{n})$. Consider $U_n=f^{-1}_n(-\frac{1}{n},\frac{1}{n})$. Assume that $U_n\neq X$. (If $f^{-1}_n(-\frac{1}{n},\frac{1}{n})=X$, then consider $U_n=B_n^{2\delta_n}$). In $\Gbs(X)$, TWO chooses $U_n$.

The play in $\Gbs(X)$ is

$\sigma(\emptyset), U_1,\dots, \sigma(U_1,\dots,U_{n-1}), U_n,\dots$

The corresponding play in $\g(\tau^s_\bk)$ is

$\psi(\emptyset), f_1,\dots, \psi(f_1,\dots,f_{n-1}),f_n, \dots$

Since $\sigma$ is a $\WS$ for ONE in $\Gbs(X)$, $\{U_n:n\in \nb\}$ is a $\gbs$-cover of $X$. It is easy to see that $\{f_n:n\in \nb\}$ converges to $\underline{0}$.
%Let $[B,\varepsilon]^s(\underline{0})$ be a neighbourhood of $\underline{0}$. Choose a $n_0$ with $\frac{1}{n_0}< \varepsilon$. There exists a sequence $\{\delta_n:n\geq n_0\}$ of positive real numbers satisfying $B^{\delta_n}\subseteq U_n=f^{-1}_n(-\frac{1}{n},\frac{1}{n})$ for all $n\geq n_0$ i.e. $f_n\in [B,\varepsilon]^s(\underline{0})$ for all $n\geq n_0$.
Hence $\psi$ is a $\WS$ for ONE in $\g(\tau^s_\bk)$.
\ep

%Note that the compactness of base in the above result is used only to prove the implication $(3)\Rightarrow (4)$.
Finally, we consider a $\PS$ for ONE. For a similar result for the space $C_p(X)$ we refer readers to \cite[Theorem 17]{ch}.

\begin{Th}
\label{Tg-3}
Let $\bk$ be a bornology with compact base on a metric space $X$. The following statements are equivalent.\\
$(1)$ $\bk$ has a countable base.\\
$(2)$ ONE has a $\PWS$ in $\g(\tau_\bk^s)$.\\
$(3)$ ONE has a $\PWS$ in $\CL(\tau^s_\bk)$.\\
$(4)$ ONE has a $\PWS$ in $\CD(\tau^s_\bk)$.\\
$(5)$ ONE has a $\PWS$ in $\obs(X)$.
\end{Th}
\bp
We prove only $(1)\Leftrightarrow (2)$.

$(1)\Rightarrow (2)$. Let $\{B_n:n\in \nb\}$ be a base for $\bk$. We assume that $B_n\subseteq B_{n+1}$ for each $n$. We define a strategy $\psi$ for ONE in $\g(\tau_\bk^s)$ as follows. For each $n$, define $\psi(n)=[B_n,\frac{1}{n}]^s(\underline{0})$. Clearly the strategy $\psi$ depends only on $n$. Now whenever $f_n\in \psi(n)$, $n\in \nb$, the sequence $\{f_n:n\in \nb\}$ converges to $\underline{0}$. Hence $\psi$ is a $\PWS$ for ONE in $\g(\tau_\bk^s)$.

$(2)\Rightarrow (1)$. Let $\psi$ be a $\PWS$ for ONE in $G(\tau_\bk^s)$. Now $\{\psi(n):n\in \nb\}$ is a sequence of open neighbourhoods of $\underline{0}$. Since $\bk$ is closed under finite unions and has closed base, we assume that $\psi(n)=[B_n,\frac{1}{n}]^s(\underline{0})$, where $B_n$ is closed for every $n\in \nb$ and $B_n\subseteq B_{n+1}$ for every $n\in \nb$. We first show that $X=\cup_{n\in \nb} B_n$. Let $x_0\in X\setminus \cup_{n\in \nb} B_n$. For each $n$ there exists a $\delta_n>0$ such that $S_{2\delta_n}(x_0)\cap B_n=\emptyset$. Choose a $f_n\in C(X)$ satisfying $f_n(x)=0$ for all $x\in \overline{B_n^{\delta_n}}$ and $f_n(x_0)=n$. Clearly for each $n$, $f_n\in [B_n,\frac{1}{n}]^s(\underline{0})$. If we choose a $B\in \bk$ with $x_0\in B$, then $f_n\not\in [B,1]^s(\underline{0})$ for any $n$. Now $\{f_n:n\in \nb\}$ is a sequence of moves by TWO which does not converge to $\underline{0}$. This contradicts that $\psi$ is a $\WS$ for ONE in $G(\tau_\bk^s)$. Therefore $X=\cup_{n\in \nb} B_n$. To complete the proof we next show that $\{B_n:n\in \nb\}$ is cofinal. Assume the contrary. Therefore there is a $B\in \bk$ satisfying $B\setminus B_n\neq \emptyset$ for each $n$. Now for each $n$ choose a $x_n\in B\setminus B_n$ and a $\delta_n>0$ with $S_{2\delta_n}(x_n)\cap B_n=\emptyset$. Again choose a $f_n\in C(X)$ satisfying $f_n(x)=0$ for all $x\in \overline{B_n^{\delta_n}}$ and $f_n(x_n)=1$. Clearly $f_n\in [B_n,\frac{1}{n}]^s(\underline{0})$ but $f_n\not\in [B,1]^s(\underline{0})$ for each $n$. Therefore $\{f_n:n\in \nb\}$ is a sequence of moves by TWO which does not converge to $\underline{0}$. Again it contradicts that $\psi$ is a $\WS$ for ONE in $G(\tau_\bk^s)$. Hence $\{B_n:n\in \nb\}$ is cofinal.
\ep

{}

\end{document}